\documentclass[12pt]{iopart}

\usepackage{iopams}  
\usepackage{graphicx}
\usepackage{color}
\usepackage{ulem}

\begin{document}

\title[Figure-eight choreographies with Lennard-Jones-type potentials]
{Figure-eight choreographies of the equal mass three-body problem with Lennard-Jones-type potentials}

\author{Hiroshi Fukuda$^1$, Toshiaki Fujiwara$^1$ and Hiroshi Ozaki$^2$}

\address{$^1$ College of Liberal Arts and Sciences, Kitasato University, 1-15-1 Kitasato, Sagamihara, Kanagawa 252-0329, Japan}
\address{$^2$ Laboratory of general education for science and technology, Faculty of Science, Tokai University, 
4-1-1 Kita-Kaname, Hiratsuka, Kanagawa, 259-1292, Japan}
\ead{fukuda@kitasato-u.ac.jp, fujiwara@kitasato-u.ac.jp and ozaki@tokai-u.jp}
\vspace{10pt}
\begin{indented}
\item[]May 2016
\end{indented}

\begin{abstract}
We report on figure-eight choreographic solutions to 
a system of three identical particles interacting through a potential of Lennard-Jones-type, 
$1/r^{12} - 1/r^6$ where $r$ is a distance between the particles. 
By numerical search, we found there are a multitude of such solutions. 
A series of them are close to a figure-eight solutions to a homogeneous system with no $1/r^{12}$ term in the potential. 
The rest are very different from them and have several points with large curvatures in their figure-eight orbits, 
at which particles are repelled.
Here figure-eight choreographies are the periodic motion 
whose shape is symmetric in both horizontal and vertical axis,  
starting with an isosceles triangle configuration and 
going back to an isosceles triangle configuration with opposite direction through Euler configuration.
Thus the lobe of this figure-eight may be complex shape and needs not to be convex.
  
\end{abstract}

%
%
%
%
%

\section{Introduction}
\label{sec:intro}
Choreographic motion of $N$ bodies is a periodic motion on a closed orbit, 
$N$ identical bodies chase each other on the orbit with equal time-spacing. 
Moore \cite{moore} found a remarkable figure-eight three-body choreographic solution 
under homogeneous interaction potential $-1/r^a$ by numerical calculations, 
where $r$ is a distance between bodies. 
Chenciner and Montgomery \cite{chenAndMont} gave a rigorous proof of its existence 
for $a=1$, i.e., for Newtonian gravity.

Sbano \cite{sbano2005}, and Sbano and Southall \cite{sbano}, 
after that, 
studied mathematically $N$-body choreographic solutions 
under 
an
inhomogeneous potential 
\begin{equation}
   u(r)=\frac{1}{r^{12}}-\frac{1}{r^6},
\label{eq:LJu}
\end{equation}
a model potential between atoms called Lennard-Jones-type potential.
For system under homogeneous potential 
as $-1/r^a$,
if there exist a periodic solution with period $T$, 
there exist scaled solutions for any period $T$. 
However, solutions for inhomogeneous potential can not be scaled.
Sbano and Southall \cite{sbano} proved that 
there exist at least two $N$-body choreographic solutions for sufficiently large period $T$, 
and there exists no solution for small period $T$. 

Choreographic three-body motion on the lemniscate \cite{lem} 
is another example for the figure-eight choreography under inhomogeneous potential
though potential 
$1/2 \ln r -(\sqrt{3}/24) r^2$ is very strange.
To our knowledge, there is no other study on the figure-eight choreography under inhomogeneous potential.   

In this paper, we study figure-eight choreographic solutions to 
a system of three identical bodies interacting through Lennard-Jones-type potential (\ref{eq:LJu})
in classical mechanics by numerical calculations.
In section \ref{sec:homo}, 
we investigate figure-eight choreographic solution under homogeneous potential $-1/r^6$,
attractive term of Lennard-Jones-type potential (\ref{eq:LJu}),
and consider a general construction of figure-eight choreography independent of the potential energy.
In section \ref{sec:inhomo},
we define figure-eight choreographic solutions under 
Lennard-Jones-type potential (\ref{eq:LJu}) 
and show a series of solutions tending toward solution to the homogeneous system.
We investigate other series of solutions 
that was predicted by Sbano and Southall
in section \ref{sec:other}.
Section \ref{sec:summary} is a summary and discussions.
Our numerical results in this paper were calculated by Mathematica 10.4 in its default precision.

\section{
Construction of figure-eight choreography
}
\label{sec:homo}
We consider solutions to a equation of motion, 
\begin{equation}
  \ddot{q}_i = -\frac{\partial U}{\partial q_i}, \; i=0,1,2,
  \label{eq:motion}
\end{equation}
for a system of three identical bodies interacting through a potential $U$
where $q_i(t)=(x_i(t), y_i(t))$ is a position vector of body $i$
in a plane of the motion, and dot represents a differentiation in $t$.

In this section, we take homogeneous potential, a power law potential as $U$, 
\begin{equation}
    U^{(a)}= -\sum_{i>j} \frac{1}{r_{ij}^a},
\label{eq:homo}
\end{equation}
where $r_{ij}=|q_i - q_j|$ is a distance between body $i$ and $j$.
For $a=1$, 
there exist figure-eight choreographic solution \cite{moore,chenAndMont},
and 
we could obtain that for $a=2,3,\ldots,14$ numerically.

\begin{figure}
   \centering
   \includegraphics[width=8cm]{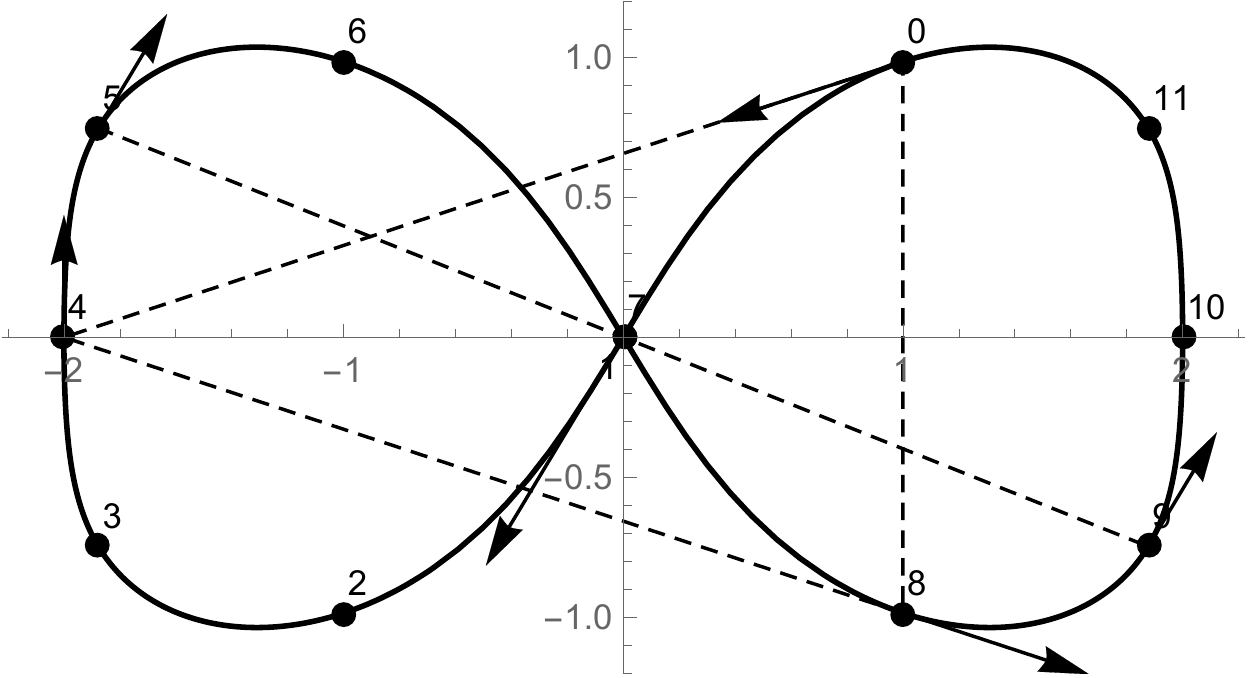} 
   \caption{
   Figure-eight choreographic solution for homogeneous potential $U^{(6)}$. 
Points labeled by ($k \bmod 12$) are the positions of bodies $q_0(k T/12)$, $q_1((k-4)T/12)$ 
and $q_2((k+4)T/12)$ where $T$ is the period of the motion and $k$ is integer. Arrows are their velocity vectors.
   }
   \label{fig:alpha6}
\end{figure}

In figure~\ref{fig:alpha6}, $q_i(t)$ of the figure-eight choreographic solution for $a=6$ is shown.
Two perpendicular mirror symmetric axes of figure-eight are taken as $x$- and $y$-axis, 
and thus the origin is the center of mass,
$\sum_i q_i =0$.
Points labeled by 
($k \bmod 12$)
are the positions $q_0(k T/12)$ where $T$ is the period of the motion 
and $k$ is integer.
Because of equal time spacing of choreography, 
\begin{equation}
  q_{i+1}(t)=q_i(t+T/3), \; i=0,1,2,
\label{eq:choreo}
\end{equation}
they are also the positions $q_1((k-4)T/12)$ and $q_2((k+4)T/12)$.
Here and hereafter it is assumed that subscripts are modulo 3.

In the time interval $T/12$, 
the
figure-eight solution takes two special configurations alternately, 
Euler configuration when one body is in the origin, and
an isosceles triangle configuration 
when one body is on the $x$-axis.
In figure~\ref{fig:alpha6}, three bodies take the isosceles triangle configuration 
at $t=2 k T/12$ and Euler configuration at $t=(2k+1) T/12$. 
The dashed triangle in figure~\ref{fig:alpha6} is the isosceles triangle at $t=0$ and
the dashed line segment the successive Euler configuration at $t=T/12$.

\subsection{
Isosceles triangle configuration
}
In the isosceles triangle configuration at $t=0$, 
we denote the position of body in the first quadrant labeled by $0$, $q_0(0)$, as $(x_0,y_0)$.
Then that in the forth quadrant labeled by $8$, $q_2(0)$, is $(x_0,-y_0)$ by the definition and 
that in the $x$-axis labeled by $4$, $q_1(0)$, is $(-2 x_0,0)$ by $\sum_i q_i=0$. 
Thus, the positions of bodies in the isosceles triangle configuration at $t=0$ are written by $(x_0, y_0)$ as
\begin{equation}
  q_0(0)=(x_0,y_0), \; q_1(0)=(-2x_0, 0), \; q_2(0)=(x_0,-y_0).
  \label{eq:isos1}
\end{equation}

In this configuration, velocity vector $\dot{q}_1(0)$,
shown by arrow at point labeled by $4$ in figure~\ref{fig:alpha6}, is parallel to $y$-axis 
since the figure-eight is symmetric in the $x$-axis and its curvature is continuous. 
By the same symmetry, 
tangent lines at $q_0(0)$ and $q_2(0)$ are symmetric in $x$-axis 
thus meet at a point $c$ in the $x$-axis. 

Here the total angular momentum $\sum_i q_i \times \dot{q}_i$ is zero 
since the areas of the left and right lobe of the figure-eight orbit are equal.
The total linear momentum $\sum_i \dot{q}_i$ is also zero since $\sum_i q_i=0$. 

By the zero total angular and linear momentum, we have
$|\sum_i q_i \times \dot{q}_i| = |\sum_i (q_i-c) \times \dot{q}_i| = |(q_1(0)-c)\times \dot{q}_1(0)|=|q_1(0)-c|\cdot|\dot{q}_1(0)|=0$.
Hence $\dot{q}_1(0)=0$ or $c=q_1(0)$. 
Suppose $\dot{q}_1(0)=0$ holds, 
either $\dot{q}_0(0)=\dot{q}_2(0)=0$ or $\dot{q}_0(0)=(0,\dot{y}_0)=-\dot{q}_2(0) \ne 0$ 
is deduced, 
which leads to a motion keeping the isosceles triangle, 
$q_0=(x(t),y(t))$, $q_1=(-2x(t),0)$ and $q_2=(x(t),-y(t))$, with $(x(0),y(0))=(x_0,y_0)$.
Therefore $c=q_1(0)$,
that is, $\dot{q}_0(0)$ and $\dot{q}_2(0)$ are parallel to the edges of isosceles triangle 
as shown in figure~\ref{fig:alpha6}.
Thus, to be $\sum_i \dot{q}_i=0$ and the motion in the left lobe clockwise,
the velocity vectors of bodies in the isosceles triangle configuration at $t=0$
are written by $(x_0, y_0)$ and $v>0$ as
\begin{equation}
\eqalign{
  \dot{q}_0(0) = \frac{v(q_1(0)-q_0(0))}{|q_1(0)-q_0(0)|} , \; 
\cr
  \dot{q}_1(0) = (0, \frac{2v}{\sqrt{1+(3x_0/y_0)^2}}), \;
\cr
  \dot{q}_2(0) = \frac{v(q_2(0)-q_1(0))}{|q_2(0)-q_1(0)|}.  
}
  \label{eq:isos2}
\end{equation}

\subsection{
Euler configuration
}
\label{sec:euler}
In the Euler configuration at $t=T/12$, 
body $0$ at $q_0(T/12)$ labeled by $1$ in figure~\ref{fig:alpha6} is in the origin.
The body $1$ at $q_1(T/12)$ labeled by $5$ is in the second quadrant and
the body $2$ at $q_2(T/12)$ labeled by $9$ in the forth quadrant
where $q_2(T/12)=-q_1(T/12)$ by $\sum_i q_i=0$.
Thus, 
the positions of bodies in the Euler configuration at $t=T/12$
are written by $q_1 \ne 0$ as
\begin{equation}
  q_0=0, \; q_2=-q_1.
  \label{eq:euler1}
\end{equation}

In this configuration, 
velocity vectors $\dot{q}_1$ and $\dot{q}_2$ are parallel
as shown by arrows at points labeled by $5$ and $9$ in figure~\ref{fig:alpha6},
since the figure-eight is symmetric in both $x$ and $y$ axis. 
Thus, writing $\dot{q}_2=a \dot{q}_1$ we have, by zero total angular momentum, 
$\sum_i q_i \times \dot{q}_i = (1-a)q_1 \times \dot{q}_1=0$ 
where $q_1 \ne 0$.
If either $\dot{q}_1=0$ or  $q_1 \parallel \dot{q}_1 \ne 0$ holds, 
one dimensional motion along a line connecting 0 and $q_1$ is deduced.
Thus $a=1$, and $\dot{q}_1=\dot{q}_2$.
Therefore, to be $\sum_i \dot{q}_i=0$, 
the velocity vectors of bodies in the Euler configuration at $t=T/12$
are written by $\dot{q}_0$ as
\begin{equation}
  \dot{q}_1 = \dot{q}_2 = -\dot{q}_0 /2.
  \label{eq:euler2}
\end{equation}
These velocity vectors are shown by arrows at points labeled by $5$, $9$ and $7$ in figure~\ref{fig:alpha6}.

\subsection{
Construction of figure-eight choreography
}
Inversely, we suppose that 
a three body motion $q_i(t)$ satisfies isosceles triangle configuration 
(\ref{eq:isos1}) and (\ref{eq:isos2})
at $t=0$,
and Euler configuration 
(\ref{eq:euler1}) and (\ref{eq:euler2})
at 
some
$t=t_0>0$.
Because the Euler conditions (\ref{eq:euler1}) and (\ref{eq:euler2}) at $t=t_0$, 
and equation of motion (\ref{eq:motion}) are invariant 
under inversions of $x$, $y$ and $t$ with exchange of bodies 1 and 2
we have 
\begin{equation}
  q_i(t+t_0)=-q_{2i}(t_0-t).
  \label{eq:sym1}
\end{equation}

By setting $t=t_0$ in the equation (\ref{eq:sym1}) and in its derivative in $t$,
we obtain $q_i(2t_0)=-q_{2i}(0)$ and $\dot{q}_i(2t_0)=\dot{q}_{2i}(0)$ at $t=2t_0$,
which are the inversions of $x$ with cyclic permutation of subscript $i$ to $i+2$ 
in the isosceles triangle conditions (\ref{eq:isos1}) and (\ref{eq:isos2}) at $t=0$.
Since the equation of motion (\ref{eq:motion}) is 
invariant under these transformations
we have
\begin{equation}
  q_i(t+2t_0)=(-x_{i+2}(t),y_{i+2}(t)).
  \label{eq:sym2}
\end{equation}

Using the equation (\ref{eq:sym2}) twice 
we obtain
the choreographic relation (\ref{eq:choreo}) with $T/3=4t_0$, 
thus the motion is choreographic and periodic with period $12t_0$.

The equation (\ref{eq:sym1}), (\ref{eq:sym2}),
\begin{equation}
  q_i(t+3t_0)=(x_{2i+1}(t_0-t), -y_{2i+1}(t_0-t))
  \label{eq:sym3}
\end{equation}
obtained by substitution of the equation (\ref{eq:sym1}) into (\ref{eq:sym2}),
and $q_i(t)$ itself
give the one third of the orbit, $q_i(t)$ for $0 \le t \le 4t_0$, 
by the four possible inversions in $x$ and/or $y$ of 
the orbit for $0 \le t \le t_0$ of the body $2i$, $i+2$, $2i+1$ and $i$.
%
The rest of the orbit, $q_i(t)$ for $4t_0 \le t \le 8t_0$ and for $8t_0 \le t \le 12t_0$ are 
given by shift of the subscripts $i$ by one and two, 
respectively,
by choreographic relation.
Each of the four subscripts, $2i$, $i+2$, $2i+1$ and $i$,
yields a set $\{0,1,2\}$ by the shifts of $i$ by zero, one and two.
Thus the orbit $q_i(t)$ for $0 \le t \le 12t_0$ is symmetric in $x$ and $y$ axis. 

Consequently three body motion $q_i(t)$ which satisfies the conditions (\ref{eq:isos1})--(\ref{eq:euler2})
are the choreography satisfying the choreographic relations (\ref{eq:choreo}) with period $T=12t_0$. 
The orbit is figure-eight shape in the sense that it is symmetric in $x$ and $y$ axis 
and passes through the origin.
This result holds for the other potential $U$ if it is invariant under inversions of $x$, $y$ and $t$.
%
Note that the set of initial conditions $q_i(0)$ and $\dot{q}_i(0)$ in equations (\ref{eq:isos1}) and (\ref{eq:isos2}) 
are determined by three parameters $(x_0,y_0,v)$. 

For the figure-eight solution under homogeneous potential (\ref{eq:homo}) with $a=6$ 
shown in figure~\ref{fig:alpha6}, the three parameters
are 
\begin{equation}
  (x_0, y_0, v )=(x_0, 0.985945 x_0, 0.234675 x_0^{-3})
\label{eq:alpha6}
\end{equation}
where $x_0>0$ and its period is $T=61.2000 x_0^4$.

\section{
Figure-eight choreography
under Lennard-Jones-type potential}
\label{sec:inhomo}
We, then consider motions under Lennard-Jones-type potential 
\begin{equation}
   U^{(12,6)}=\sum_{i>j} u(r_{ij}),
\label{eq:LJ}
\end{equation}
numerically.
We define the figure-eight choreography 
as a motion starting with the isosceles triangle configurations 
(\ref{eq:isos1}) and (\ref{eq:isos2}),
and going through the Euler configuration 
(\ref{eq:euler1}) and (\ref{eq:euler2}).

We start numerical integration at $t=0$ 
with the isosceles triangle initial conditions 
(\ref{eq:isos1}) and (\ref{eq:isos2})
determined by a set of parameters $(x_0, y_0, v)$.
We stop integration if three bodies are aligned in a line, that is,
\begin{equation}
   (q_1-q_0) \times (q_2-q_0) = 0
\label{eq:collinear}
\end{equation}
is satisfied.
We denote this instant as $t_f$. 
%
%
%
Then we investigate the outer product
\begin{equation}
   P(x_0,y_0,v)=\dot{q_2} \times \dot{q_1},
\end{equation}
and the difference 
\begin{equation}
  D(x_0,y_0,v)=(q_1-q_0)^2- (q_2-q_0)^2, 
\end{equation}
at $t=t_{f}$.

Thus the Euler configuration (\ref{eq:euler1}) and (\ref{eq:euler2}) are written by
$P(x_0,y_0,v)=0$ and $D(x_0,y_0,v)=0$ 
since $D(x_0,y_0,v)=0$ with (\ref{eq:collinear}) leads to $q_0(t_f)=0$, and 
$P(x_0,y_0,v)=0$ 
with the zero total angular momentum 
assured by the isosceles triangle conditions (\ref{eq:isos1}) and (\ref{eq:isos2})
leads to
$\dot{q}_1(t_f) = \dot{q}_2(t_f)$
according to the argument in the section \ref{sec:euler}.
Therefore if these conditions are satisfied 
the solution has choreographic properties  (\ref{eq:choreo}) with period $T=12t_{f}$
and has the symmetry of figure-eight, 
as stated in section \ref{sec:homo}. 

\begin{figure}
   \centering
   \includegraphics[width=9cm]{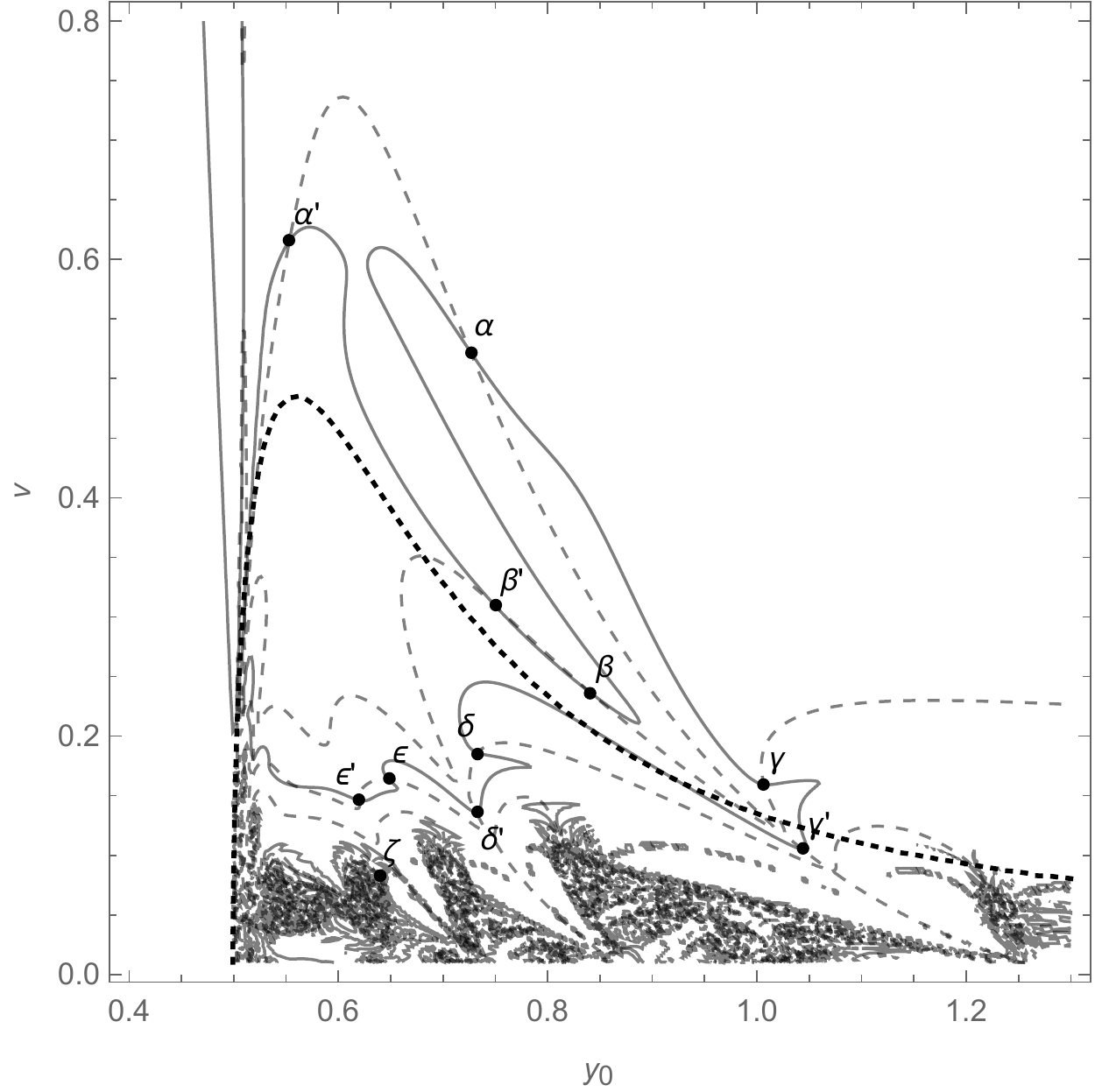} 
   \caption{
   $P(x_0,y_0,v)$=0 curves and $D(x_0,y_0,v)=0$ curves in the $(y_0,v)$-plane for $x_0=0.75$. 
   Solid curves are $D(x_0,y_0,v)=0$ and dashed curves $P(x_0,y_0,v)=0$.
   Dotted curve is $E(x_0,y_0,v)=0$.
   }
\label{fig:contour}
\end{figure}

In figure~\ref{fig:contour},
for $x_0=0.75$,
$P(x_0,y_0,v)=0$ and $D(x_0,y_0,v)=0$ curves in $(y_0,v)$-plane are shown
by dashed and solid curves, respectively.
The crossing points of the dashed and solid curves are 
the parameters $(y_0,v)$ for figure-eight choreographic solutions. 

We discuss the solution $(x_0, y_0,v)=(0.75, 0.725966,0.522742)$ labeled by $\alpha$ in figure~\ref{fig:contour}.
This solution is very close to the solution to the homogeneous system (\ref{eq:alpha6})
shown in figure~\ref{fig:alpha6}.
If we plot orbits of two solutions in the same figure, 
it is difficult to distinguish them in the usual printing resolution.

\begin{figure}
   \centering
   \includegraphics[width=7cm]{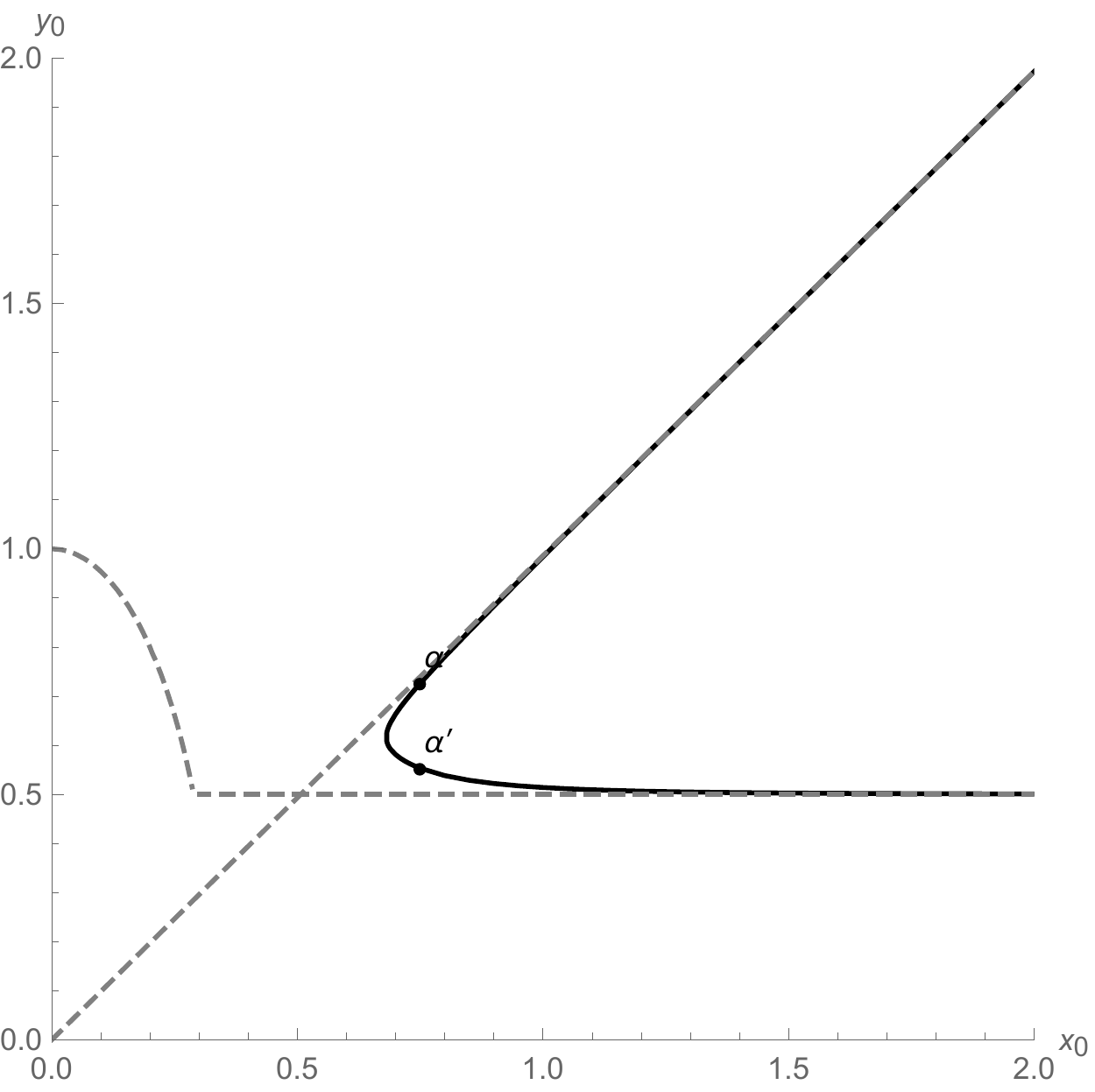} 
   \caption{
A series of solution $\alpha$.
Solid curve is a set of $(x_0, y_0)$ for the series of solution $\alpha$.
Points labeled by $\alpha$ and $\alpha'$ are the solutions labeled by the same labels in figure~\ref{fig:contour}, respectively.
The dashed line is the solution to the homogeneous system (\ref{eq:homoline}). 
The other dashed curve is a boundary of the region of strong repulsive force (\ref{eq:strong}).
   }
   \label{fig:aa}
\end{figure}

We follow this solution 
changing $x_0$ moderately. 
Then, we obtain a series of solution. 
In figure~\ref{fig:aa}, a set of $(x_0, y_0)$ for the series of this solution $\alpha$ is shown by solid curve.
Point labeled by $\alpha$ is the solution $\alpha$.
The dashed line, 
\begin{equation}
  y_0=0.985945 x_0 
\label{eq:homoline}
\end{equation}
is the solution to the homogeneous system (\ref{eq:alpha6}). 
For $x_0 \ge 0.75$, the solid curve 
is almost on this line as noted above.

\begin{figure}
   \centering
   \includegraphics[width=8cm]{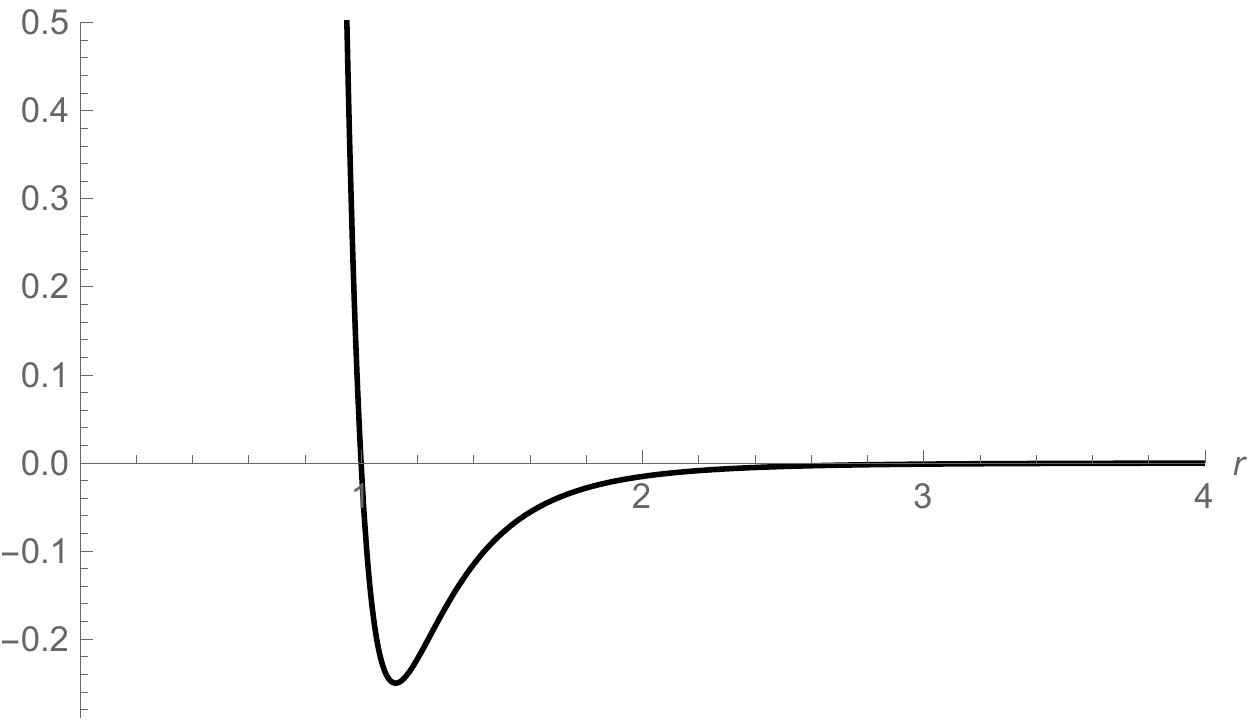} 
   \caption{
Lennard-Jones-type potential (\ref{eq:LJu}) 
has minimum at $r=r_0=2^{1/6}=1.12$ and is zero at $r=r_0'=1$.
   }
   \label{fig:LJ}
\end{figure}

The other dashed curve in figure~\ref{fig:aa} 
is a boundary of a region of strong repulsive force,
\begin{equation}
  \min(2 y_0, \sqrt{9 x_0 ^2 + y_0^2}) \le r_0' = 1,
\label{eq:strong}
\end{equation} 
in which two bodies in the isosceles triangle configuration feel strong repulsive force,
where 
$r_0'$ is a distance $r$ which makes the potential (\ref{eq:LJu}) zero.
See figure~\ref{fig:LJ}.
%
Deep in the region (\ref{eq:strong}) no solution is expected to exist 
since it is difficult to take the isosceles triangle configuration against the strong repulsive force.
Actually, in figure~\ref{fig:contour} 
no crossing point can be seen deep in the region (\ref{eq:strong}),  that is, $y_0<r_0'/2=0.5$.

Dashed line and dashed curve intersect at $x_0=0.507$.
Thus this series of solution can not exist for $x_0 < 0.507$,
and the solid curve
turns and goes along the boundary of the region (\ref{eq:strong}). 
The smallest $x_0$ for this series of solution is about $x_0=0.6812$ with $y_0=0.617578$.

\begin{figure}
   \centering
   \includegraphics[width=8cm]{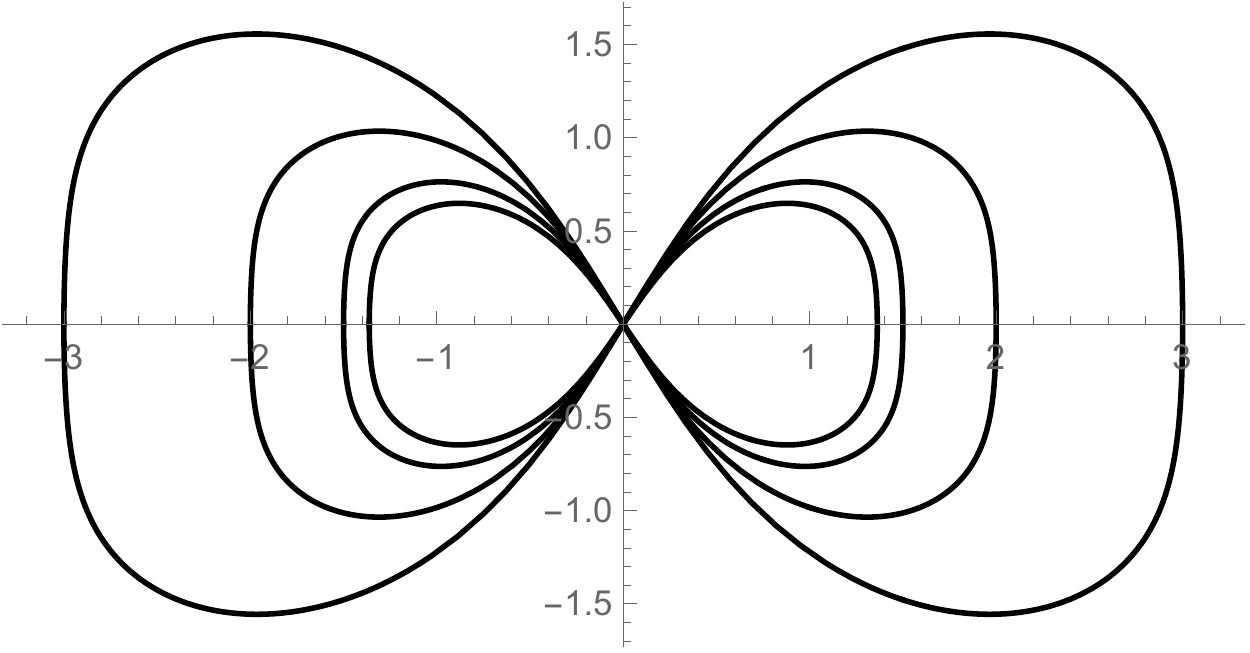} 
   \caption{
Orbits for the solutions with 
$(x_0,y_0,v)=$ 
$(0.6812, 0.617578, 0.660119)$, 
$(0.75, 0.725966, 0.522742)$,
$(1.0, 0.983588, 0.232296)$,
$(1.5, 1.478621, 0.0694721)$
in the series of solution $\alpha$. 
   }
   \label{fig:a}
\end{figure}

\begin{figure}
   \centering
   \includegraphics[width=8cm]{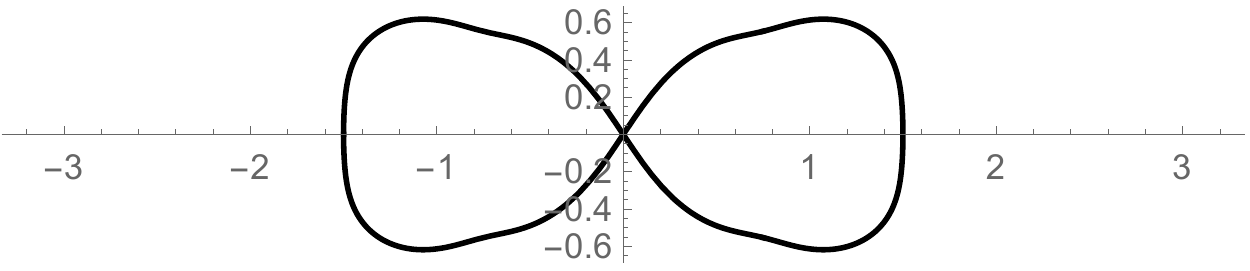} 
   \includegraphics[width=8cm]{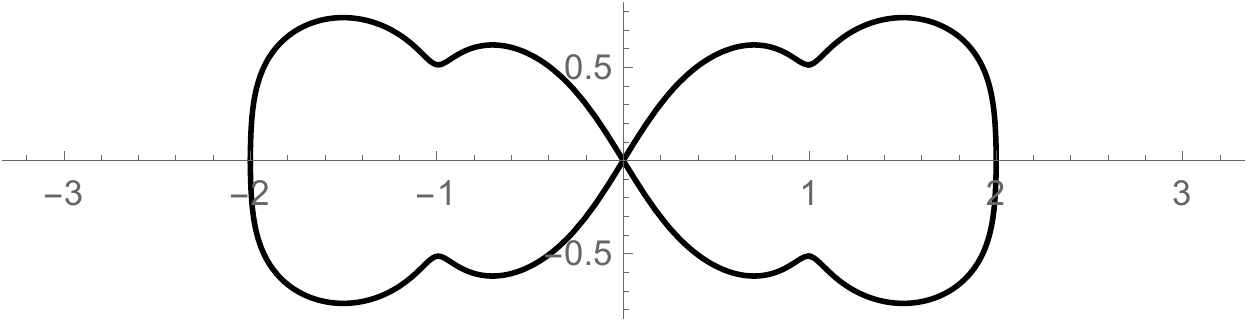} 
   \includegraphics[width=8cm]{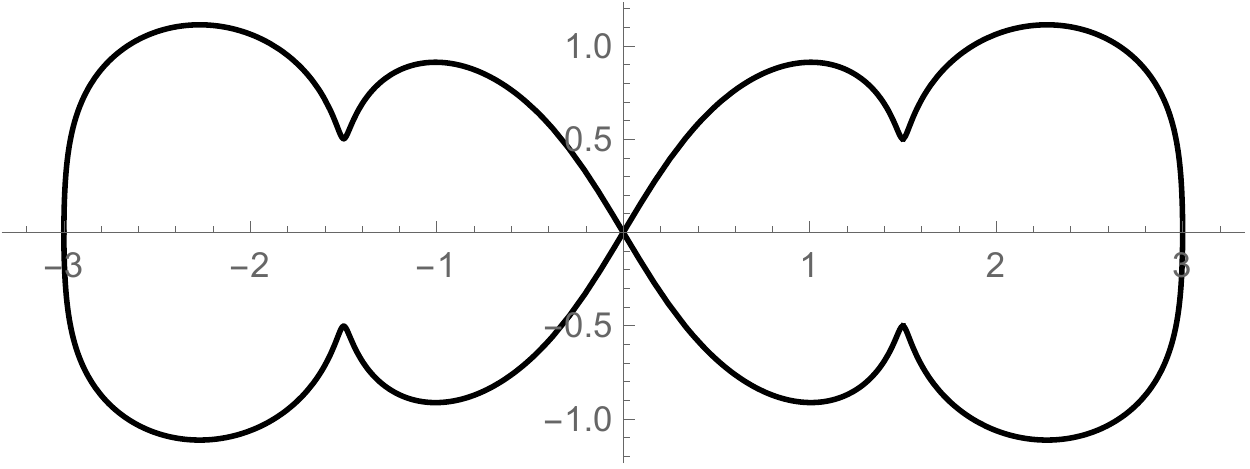} 
   \caption{
Orbits for the solutions with 
$(x_0,y_0,v)=$ 
$(0.75, 0.553223, 0.615805)$,
$(1.0, 0.513969, 0.396537)$,
$(1.5, 0.502649, 0.181062)$
in the series of solution $\alpha$. 
   }
   \label{fig:ap}
\end{figure}

In figure~\ref{fig:a}, 
orbits $q_i(t)$'s in the $y_0 \ge 0.617578$ branch of this series are shown 
for $x_0=$ 0.6812, 0.75, 1.0, 1.5.
%
They are almost scalable.
In figure~\ref{fig:ap}, those in the $y_0 < 0.617578$ branch are shown.
%
They are gourd-shaped figure-eight with necks 
slightly inside of
$x= \pm x_0$ of width about $2 y_0$.
At the both sides of 
the isosceles triangle configuration, 
two bodies are repelled since $2 y_0 < r_0$,
where $r_0=2^{1/6}=1.12$ is the distance $r$ 
which makes the potential (\ref{eq:LJu}) minimum.
See figure~\ref{fig:LJ}.   

In default precision calculation by Mathematica 
we can not get this gourd-shaped figure-eight for $x_0>2.5$.
This limit is probably due to the inaccuracy in numerical integration and
it will exist for any large $x_0$.

In figure~\ref{fig:contour}, solutions to be discussed are labeled by 
the Greek alphabets or those with prime,
where the same letters mean that the solutions belong to the same series of solution.
The solution $\alpha'$ in figure~\ref{fig:contour}
whose orbit is shown at the top in figure~\ref{fig:ap},
thus, belongs to the series of solution $\alpha$, as shown by the point $\alpha'$ in figure~\ref{fig:aa}.

\section{Other 
figure-eight choreography
}
\label{sec:other}
\begin{figure}
   \centering
   \includegraphics[width=8cm]{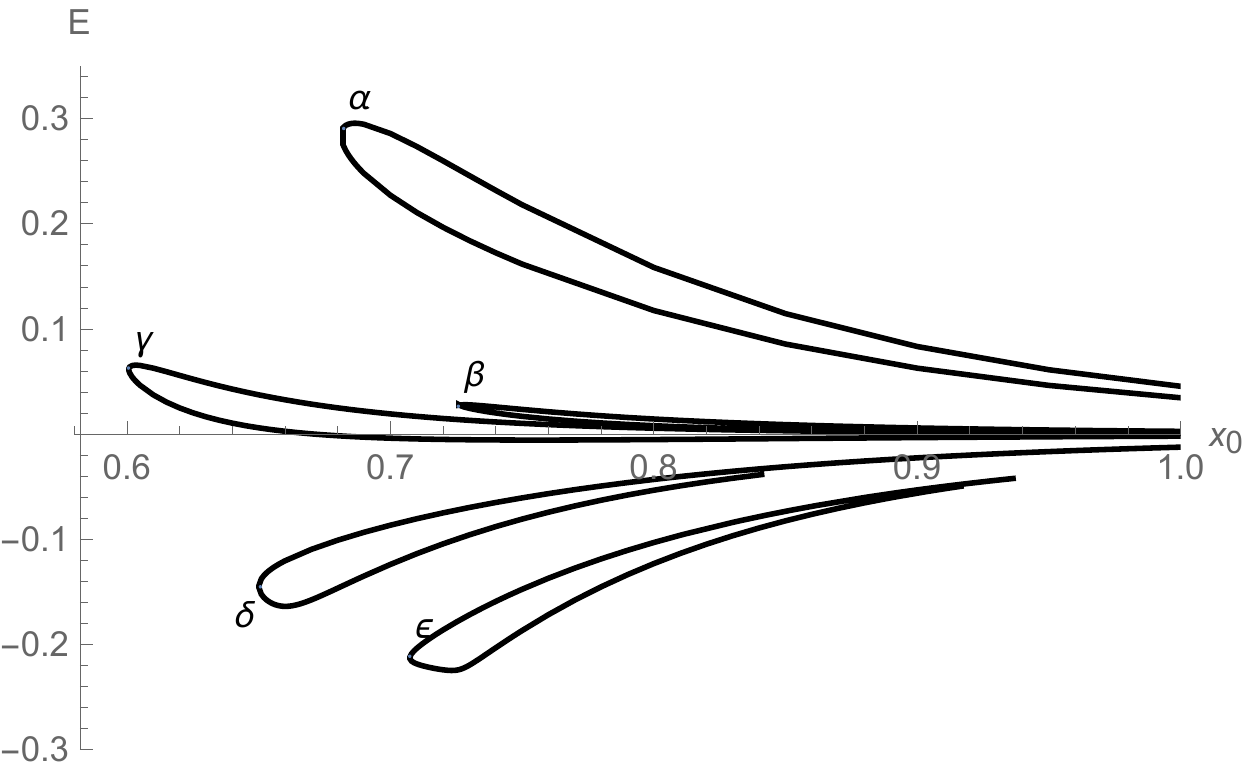} 
   \caption{
  The set $(x_0,E)$ of the solutions $\alpha$ -- $\epsilon$ in the $(x_0,E)$-plane.
   }
   \label{fig:ex0}
\end{figure}

We investigate the other series of solutions labeled by $\beta$ -- $\epsilon$ in figure~\ref{fig:contour}.
To distinguish these series we adopt total energy
\begin{equation}
  E(x_0,y_0,v)=\sum_{i} \frac{1}{2} \dot{q}_i^2+U
\label{eq:E}
\end{equation}
instead of $y_0$.
In figure~\ref{fig:ex0}, sets of $(x_0, E)$ for the series of solutions $\alpha$ -- $\epsilon$ are shown in $(x_0, E)$-plane.

In figures~\ref{fig:b}, \ref{fig:c}, \ref{fig:d} and \ref{fig:e}, 
orbits for the series of solutions $\beta$, $\gamma$, $\delta$ and $\epsilon$ are shown, respectively.
In these figures, three orbits are displayed as follows:
Orbit for the smallest $x_0$ is placed in figure (a).
For large $x_0$ orbit in higher $E$ branch is placed in figure (b) and 
that in lower $E$ branch in figure (c).

\subsection{Points of large curvature in orbits}
\begin{figure}
  \centering
  \begin{minipage}{6.5cm}
   \centering
   \includegraphics[width=6cm]{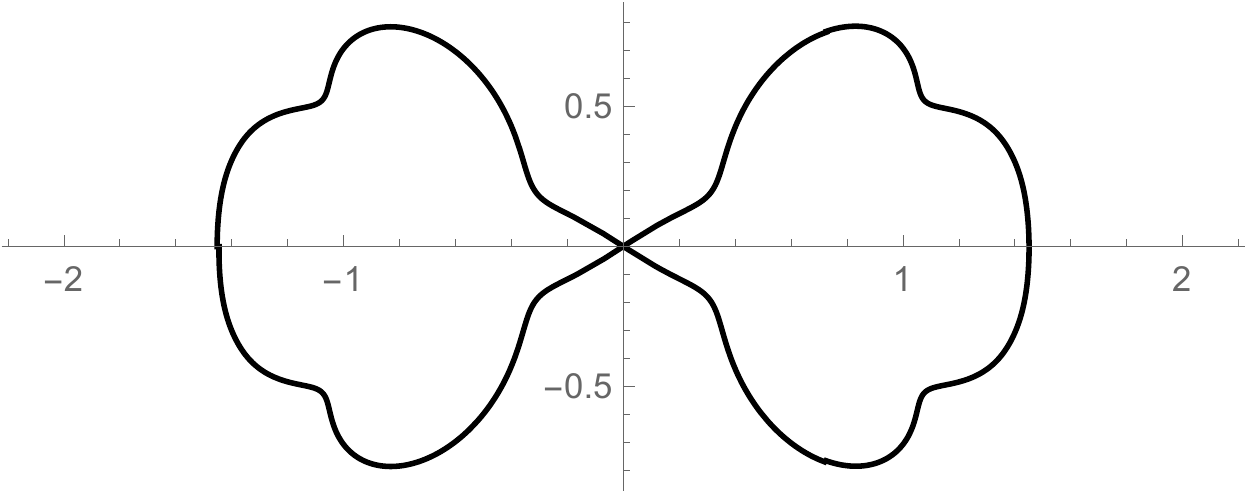}
   \scriptsize
   {\bf (a)} $(0.726, 0.766265, 0.302694; 0.0274632)$ 
  \end{minipage}
  \\
  \begin{minipage}{12.5cm}
   \centering
   \includegraphics[width=6cm]{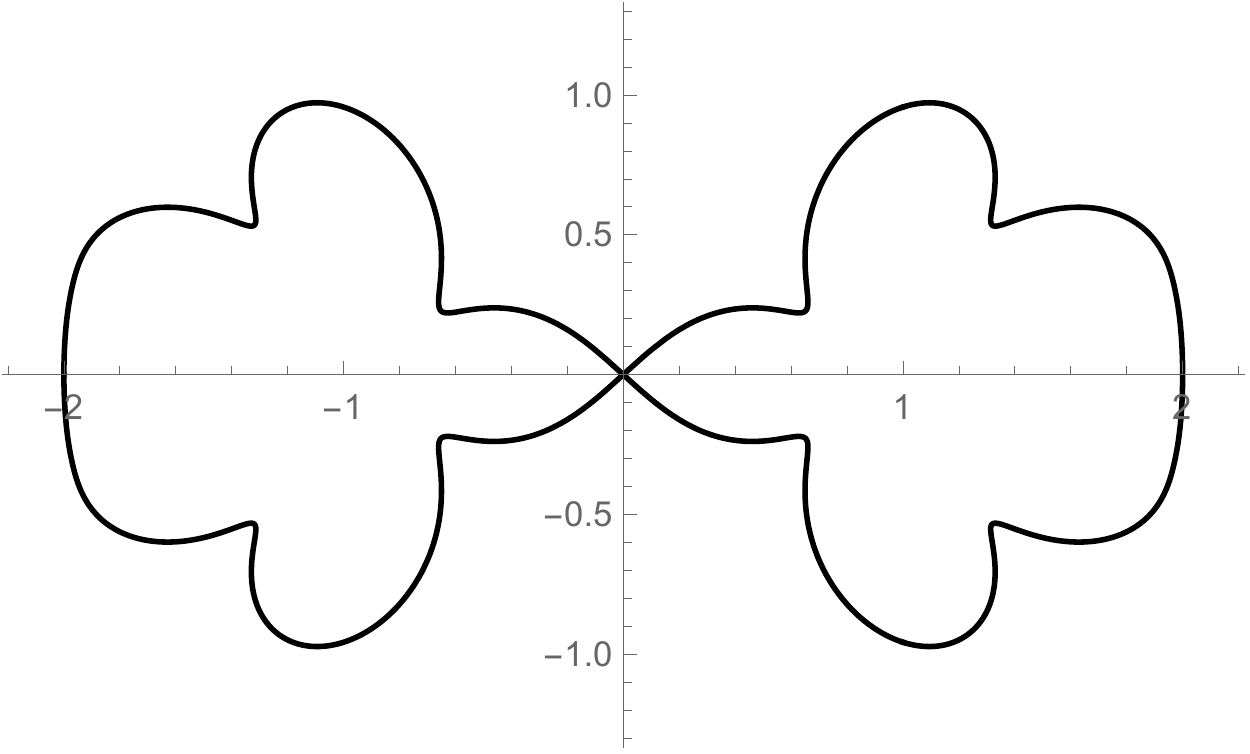} 
   \includegraphics[width=6cm]{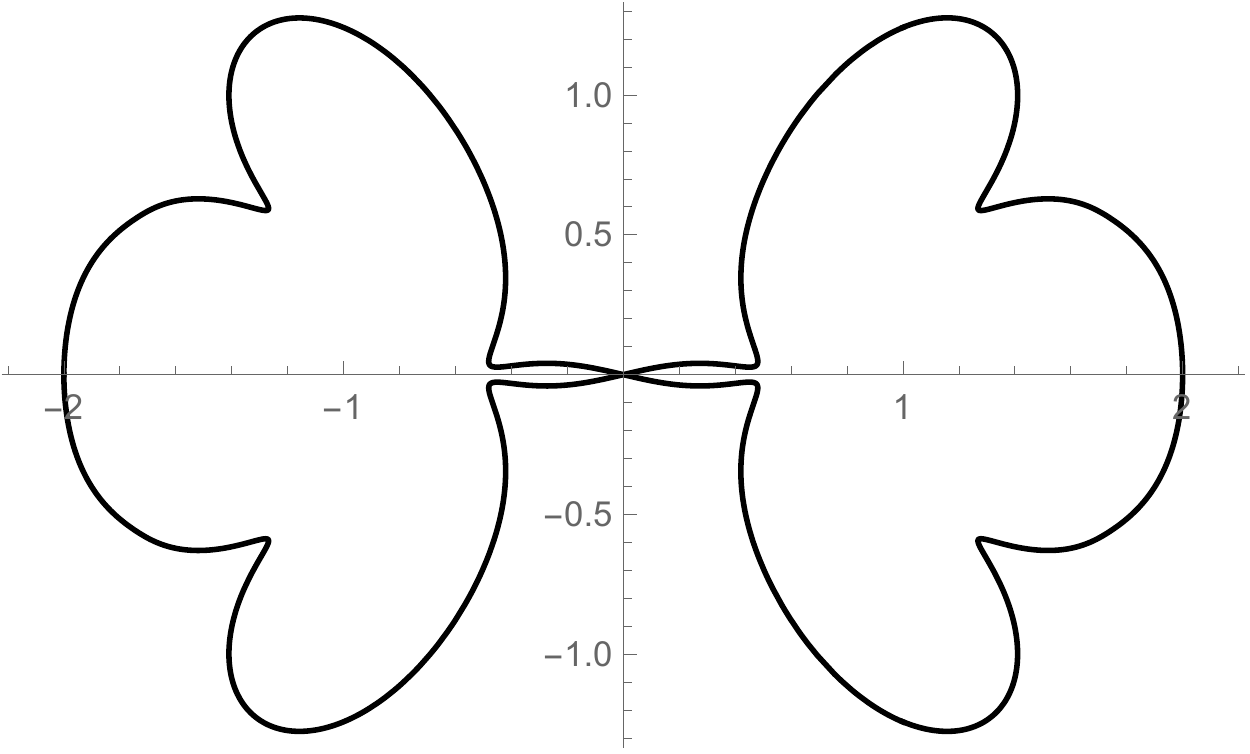} 
   \\
   \scriptsize
   {\bf (b)}  $(1.0, 0.956733, 0.144241; 0.00263843)$ \hspace{2em} 
   {\bf (c)}  $(1.0, 1.241130, 0.0717890; 0.000697053)$
  \end{minipage}
   \caption{
Orbits for a series of solution $\beta$. 
Figures in parentheses are $(x_0, y_0, v; E)$.  
   }
   \label{fig:b}
\end{figure}

\begin{figure}
  \centering
  \begin{minipage}{6.5cm}
   \centering
   \includegraphics[width=6cm]{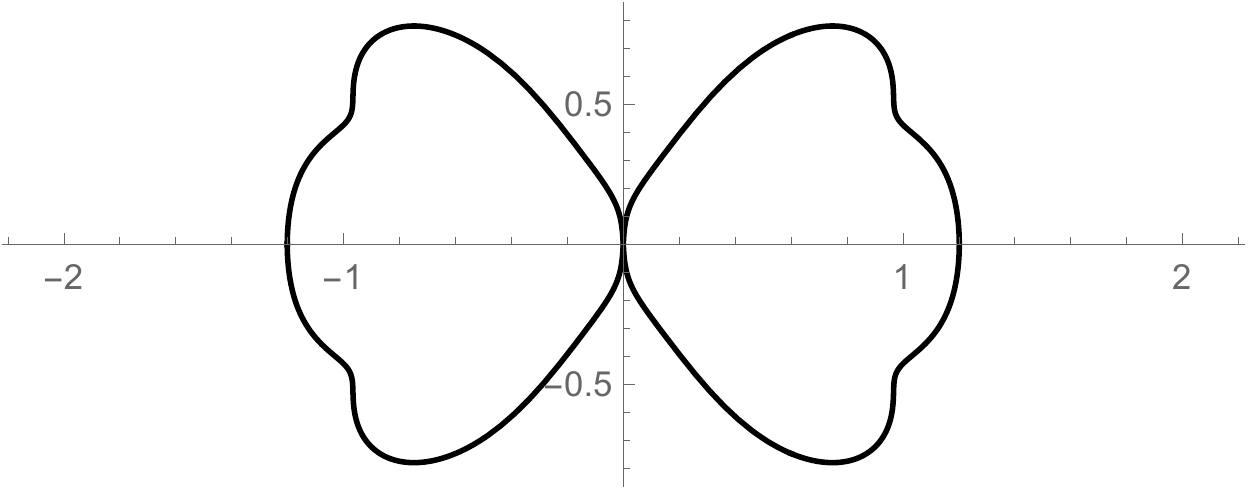}
   \scriptsize
   {\bf (a)} $(0.6007, 0.748371, 0.371779; 0.0622734)$
  \end{minipage}
  \\
  \begin{minipage}{12.5cm}
   \centering
   \includegraphics[width=6cm]{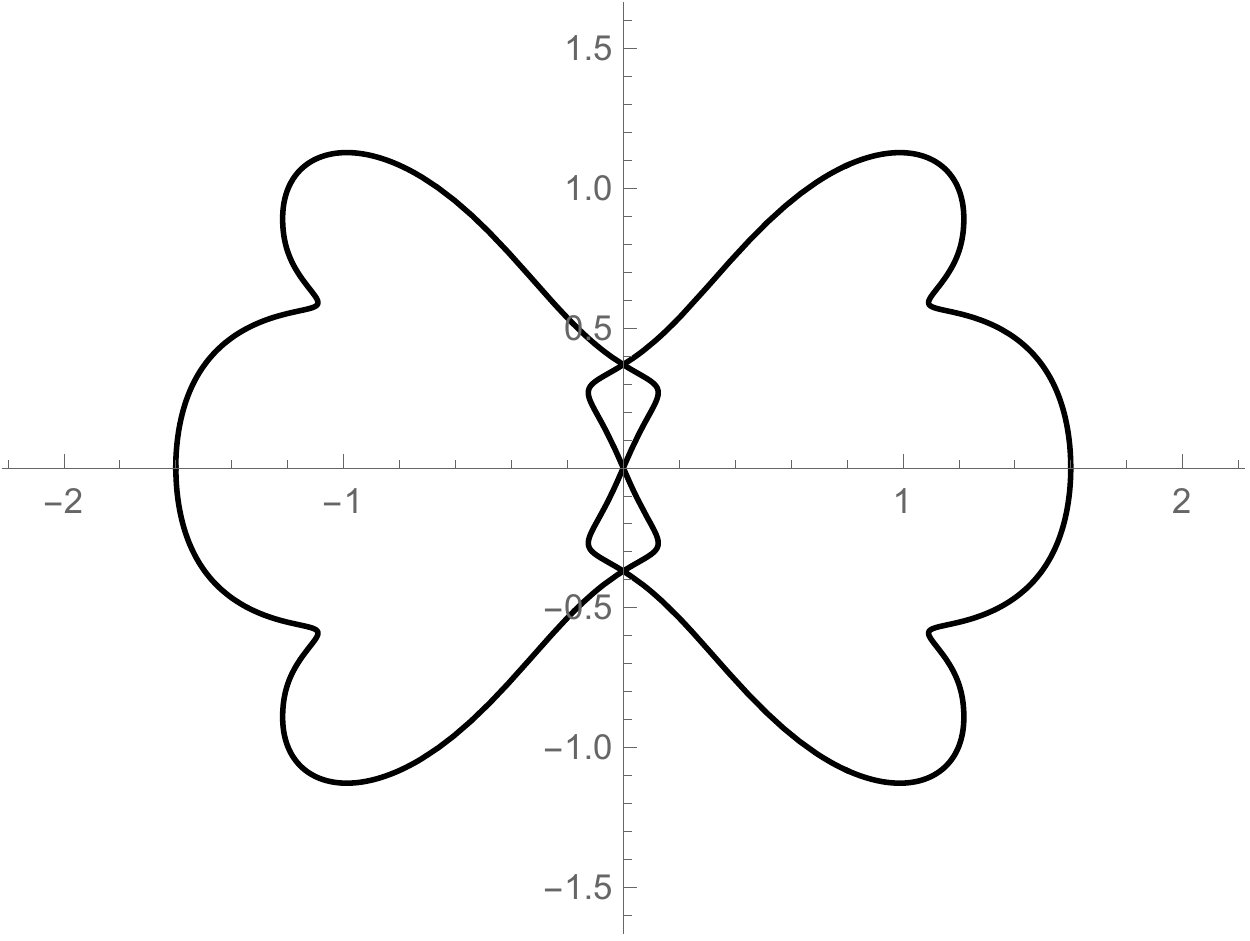} 
   \includegraphics[width=6cm]{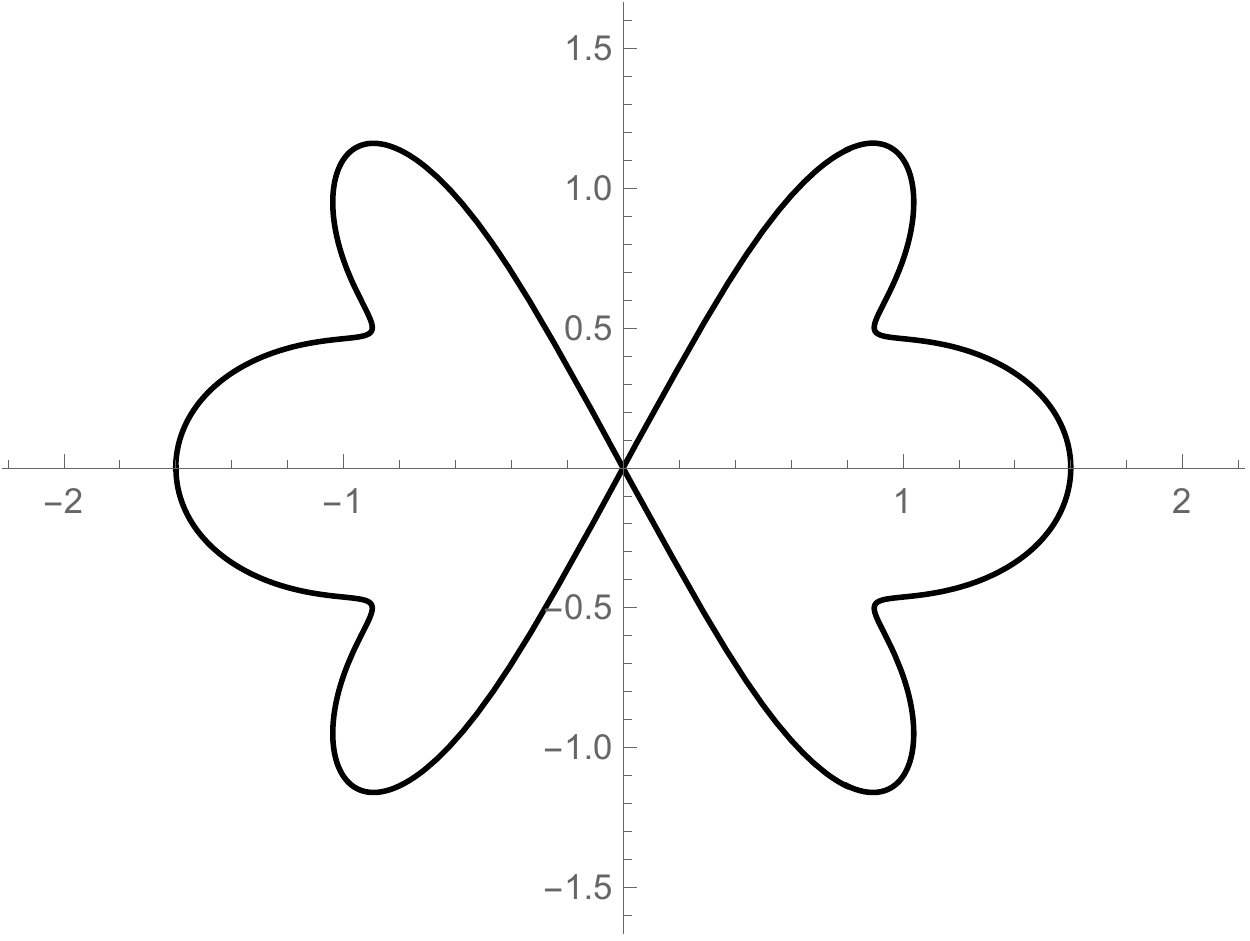} 
   \\
   \scriptsize
   {\bf (b)} $(0.8, 1.081836, 0.126051; 0.00561328)$ \hspace{2em}
   {\bf (c)} $(0.8, 1.136739, 0.0749665; -0.00519619)$
  \end{minipage}
   \caption{
Orbits for a series of solution $\gamma$. Figures in parentheses are $(x_0, y_0, v; E)$.  
   }
   \label{fig:c}
\end{figure}

\begin{figure}
  \centering
  \begin{minipage}{6.5cm}
   \centering
   \includegraphics[width=6cm]{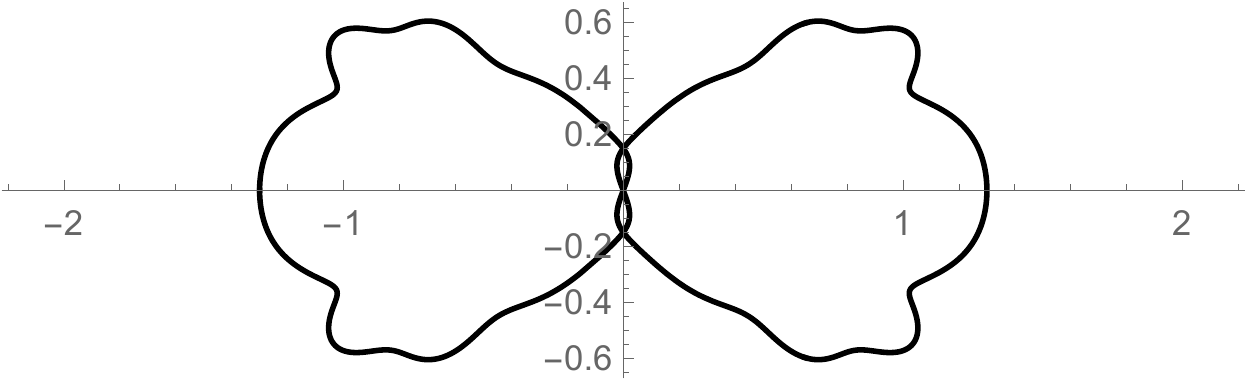}
   \scriptsize
   {\bf (a)} $(0.6501, 0.597985, 0.304229; -0.143858)$
  \end{minipage}
  \\
  \begin{minipage}{12.5cm}
   \centering
   \includegraphics[width=6cm]{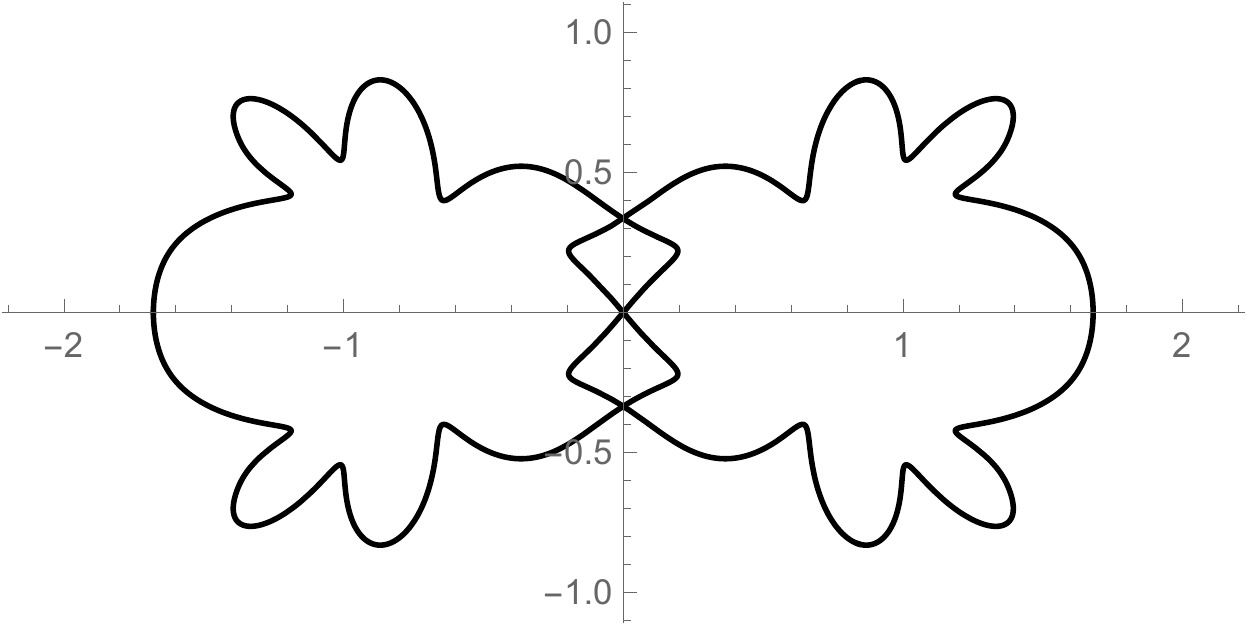} 
   \includegraphics[width=6cm]{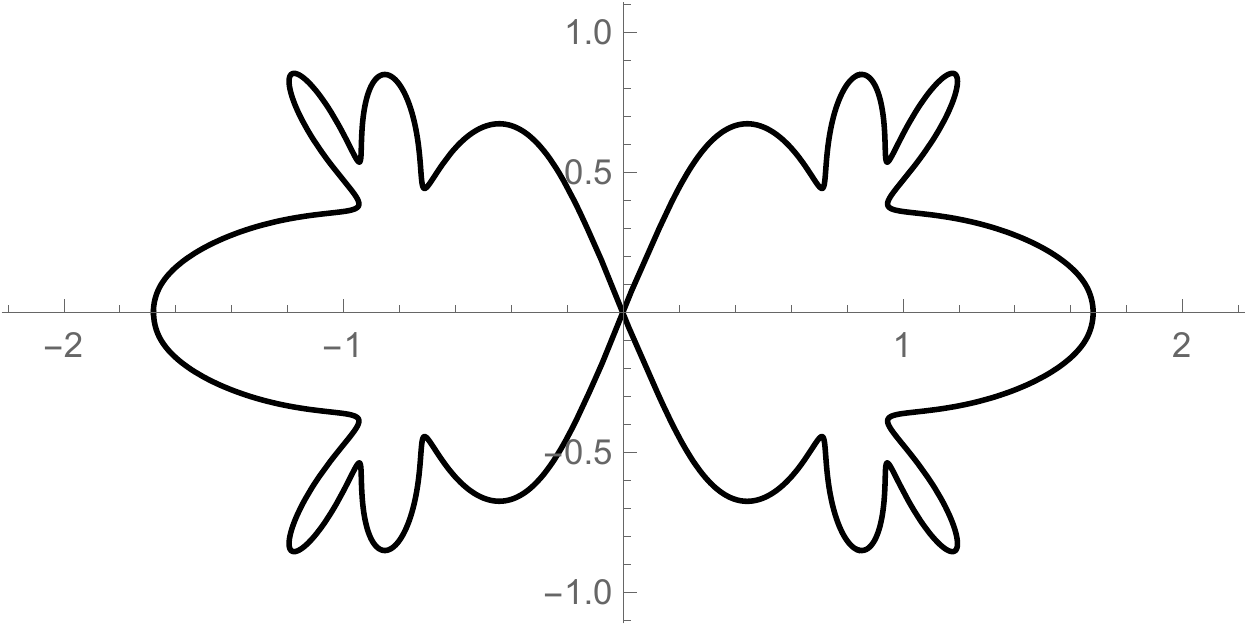} 
   \\
   \scriptsize
   {\bf (b)} $(0.84, 0.827038, 0.126408; -0.0330865)$ \hspace{2em}
   {\bf (c)} $(0.84, 0.848830, 0.0757119; -0.0387688)$
  \end{minipage}
   \caption{
Orbits for a series of solution $\delta$. Figures in parentheses are $(x_0,y_0,v; E)$.   
   }
   \label{fig:d}
\end{figure}

\begin{figure}
   \centering
  \begin{minipage}{6.5cm}
   \centering
   \includegraphics[width=6cm]{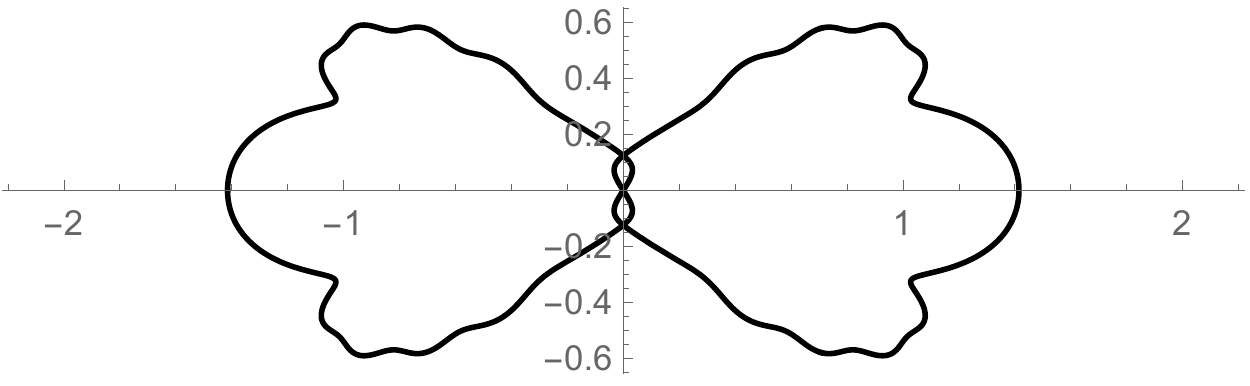}
   \scriptsize
   {\bf (a)} $(0.7074, 0.579781, 0.204620; -0.211945)$ 
  \end{minipage}
  \\
  \begin{minipage}{12.5cm}
   \centering
   \includegraphics[width=6cm]{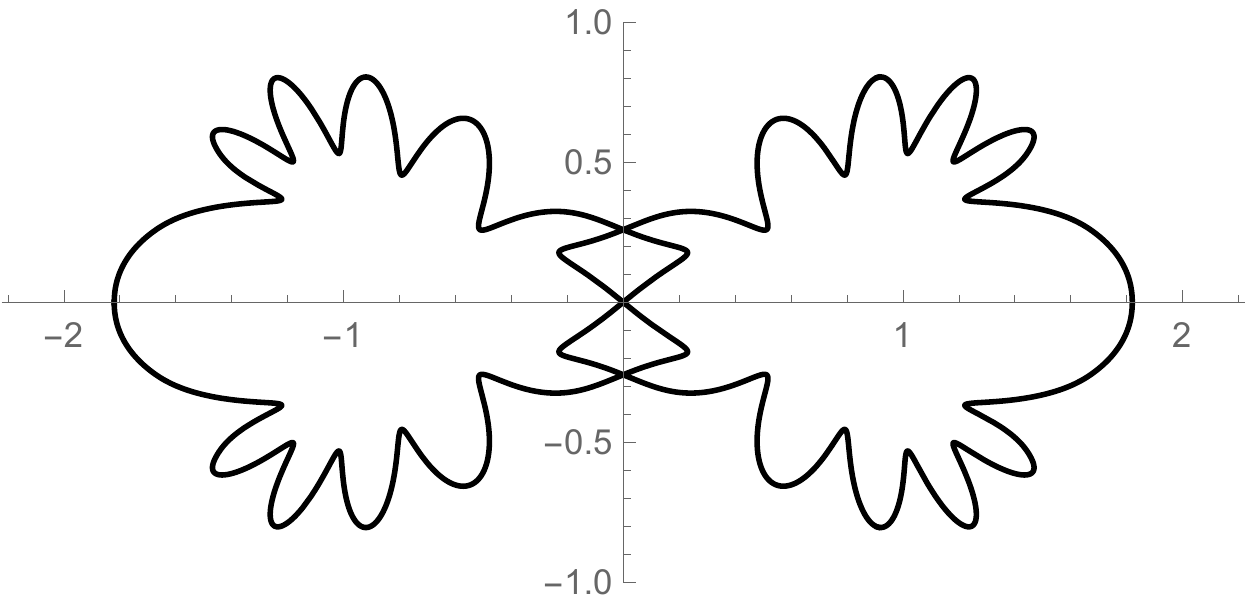} 
   \includegraphics[width=6cm]{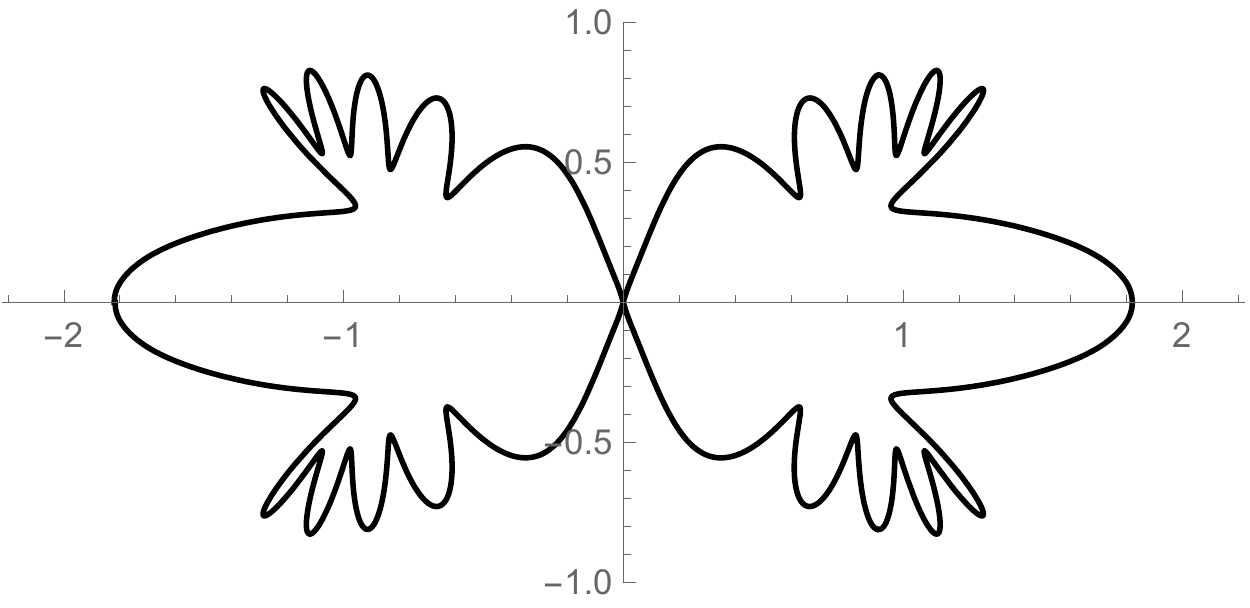} 
   \\ 
   \scriptsize 
   {\bf (b)} $(0.91, 0.803912, 0.0857343; -0.0497687)$ \hspace{2em}
   {\bf (c)} $(0.91, 0.811359, 0.0540501; -0.0521151)$
  \end{minipage}
   \caption{
   Orbits for a series of solution $\epsilon$.  Figures in parentheses are $(x_0,y_0,v; E)$.  
   }
   \label{fig:e}
\end{figure}

We notice points of large curvature in the orbits shown 
in figures~\ref{fig:ap}, \ref{fig:b}, \ref{fig:c}, \ref{fig:d} and \ref{fig:e}. 
The large curvature occurs if distance between two bodies is less than $r_0$.
We call a closed interval in $t$ which satisfies $r_{ij}(t) \le r_0$ as a collisional interval between body $i$ and $j$.
We then define number of collisions 
$n_{ij}$ by the number of collisional intervals 
between body $i$ and $j$ in $0 \le t <T$ with periodic boundary conditions. 
The number of collisions body $i$ experiences is 
\begin{equation}
  n_{i}=\sum_{j \ne i} n_{ij}.
\label{eq:ni}
\end{equation}
Since $n_0=n_1=n_2$, we use $n_0$ as the number of collisions.
We note that the term `collision' here does not mean usual collision of two point bodies $i$ and $j$ without size
in classical mechanics at which $r_{ij}=0$
but a repulsion which is considered as collision of two balls with radius $r_0/2$ 
where $r_{ij}$ can not be zero by the strong repulsive force. 

\begin{figure}
   \centering
   \includegraphics[width=8cm]{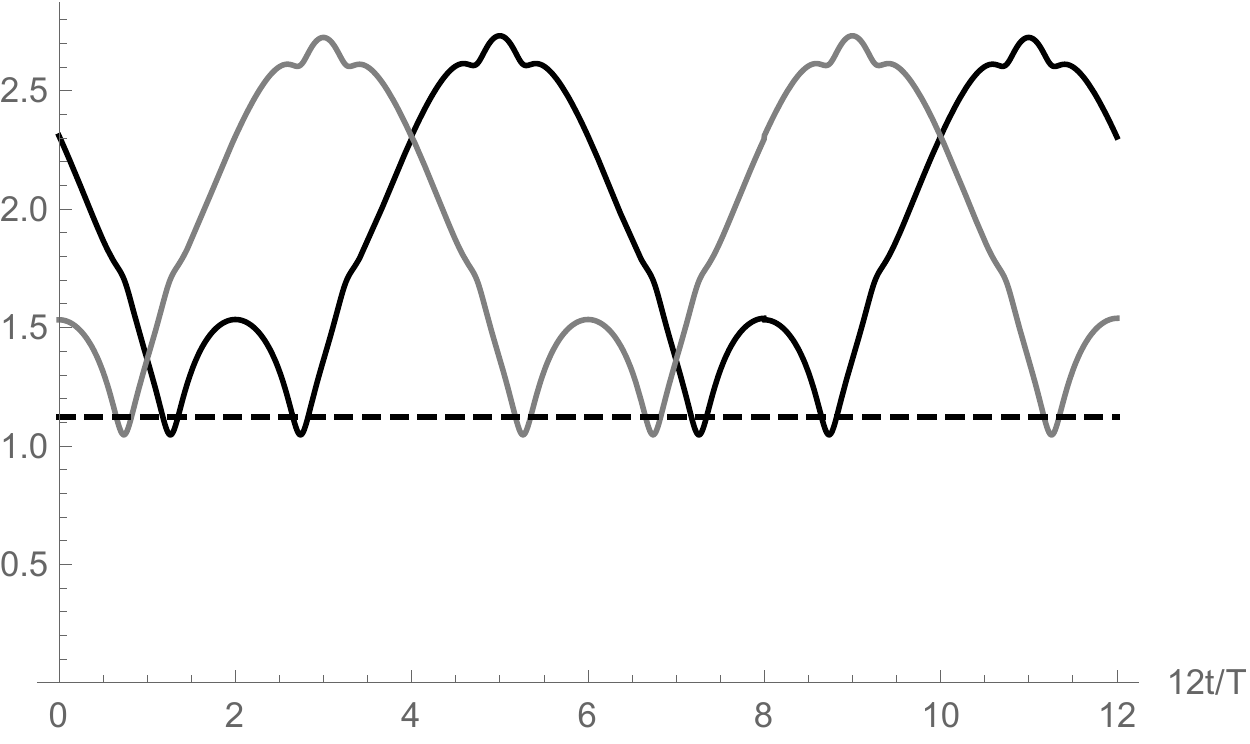} 
   \caption{
  Distances between bodies for the orbit shown in figure~\ref{fig:b} (a). 
  Black curve is $r_{01}(t)$ and gray curve is $r_{02}(t)$.
  Horizontal axis is $12t/T$. Dashed line shows $r_0$.
   }
   \label{fig:b1d}
\end{figure}

In figure~\ref{fig:b1d}, $r_{01}(t)$ and $r_{02}(t)$ for the orbit in figure~\ref{fig:b} (a) are shown.
We can count $n_{01}=4$ and $n_{02}=4$ then $n_0=n_{01}+n_{02}=8$. 
%
Similarly we can count 
$n_0$ for the other orbits.
They seem to be conserved within the same series of solutions. 
The count $n_0$ for series of solutions $\beta$ -- $\epsilon$ are 
8, 8, 16 and 24, respectively.

However, for the series of solution $\alpha$ the conservation of $n_0$ does not hold.  
The count $n_0$ for the orbits shown in figure~\ref{fig:ap} is $4$, 
and $n_0$ for those in figure~\ref{fig:a} is $0$, and $n_0$ changes around $x_0=0.72$.
Instead of the $n_0$, however, 
a number of local minimums in $r_{01}(t)$ and $r_{02}(t)$ for the orbits in figure~\ref{fig:a} is 4,
which are lowered and change to four collisional intervals.

\begin{figure}
   \centering
   \includegraphics[width=8cm]{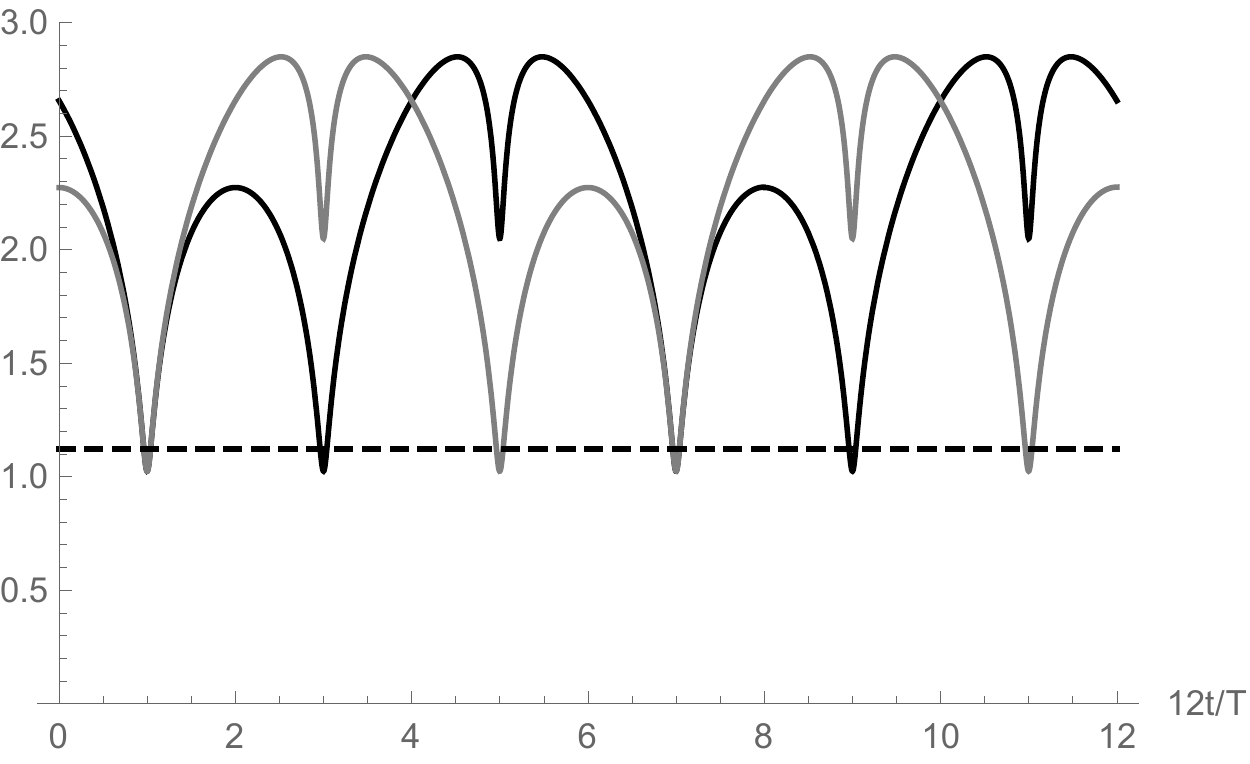} 
   \caption{
  Distances between bodies for the orbit shown in figure~\ref{fig:c} (c). 
  Black curve is $r_{01}(t)$ and gray curve is $r_{02}(t)$.
  Horizontal axis is $12t/T$. Dashed line shows $r_0$.
   }
   \label{fig:cpd}
\end{figure}

Two collisions may occur simultaneously.
In figure~\ref{fig:cpd}, 
$r_{01}(t)$ and $r_{02}(t)$ for the orbit in figure~\ref{fig:c} (c) are shown.
The collisional intervals in $r_{01}(t)$ and $r_{02}(t)$ 
are overlapped  at $t=T/12$ and $7T/12$ and 
represent simultaneous collisions between body 0 and 1, and between 0 and 2 in the Euler configuration,
where body 0 is at the origin and is collided by bodies 1 and 2 from the opposite side.

If there is no simultaneous collision,
the collisional segment, 
which we define here as a segment of the orbit corresponding to collision interval in $t$,
includes the points with large curvatures
whose accelerations are toward the outside of the lobe of figure-eight.
For simultaneous collision,
as we can see from the orbit near the origin in figure~\ref{fig:c} (c), 
the collisional segment does not include the point of large curvature and looks smooth. 

Except for the series of the solutions $\alpha$ and $\beta$ 
all series of the solutions $\gamma$ -- $\epsilon$ have simultaneous collisions.
At large $x_0$ in higher $E$ branches, 
the solutions shown in figures~\ref{fig:c} (b), \ref{fig:d} (b) and \ref{fig:e} (b),
have four collisional segments surrounding the origin. 
At the smallest $x_0$,
solutions shown in the same figures (a)
still have the four separate collisional segments. 
At large $x_0$ in lower $E$ branches, 
the four collisional segments are merged into two smooth collisional segments at the origin
shown in the same figures (c).
Though the collisional segments at the origin
look smooth,
minimum of $r_{01}(t)$ and $r_{02}(t)$ do not occur simultaneously and 
they still at slightly different time. 
They seem to occur simultaneously 
at $x_0 \to \infty$.


\subsection{Behaviors and range of the total energy}

There are three series of solutions $\alpha$, $\beta$ and $\gamma$ in $E \ge 0$ region in figure~\ref{fig:ex0}.
In figure~\ref{fig:contour}, $E=0$ curve is added by dotted curve, and 
we can see there are five solutions $\alpha$, $\alpha'$, $\beta$, $\beta'$, and $\gamma$ 
above the $E=0$ curve, i.e., $E \ge 0$ region. 
We explored such map as figure~\ref{fig:contour} for various $x_0$ with various range of $y_0$ and $v$,
and concluded numerically that 
these series of solutions $\alpha$, $\beta$ and $\gamma$ are the only series of solutions in $E \ge 0$ region.

The whole series of solutions $\alpha$ and $\beta$ exist in the region $E>0$, 
whereas the lower energy branch for the series of solution $\gamma$ 
crosses the $E=0$ once at $x_0=0.671188$ in figure~\ref{fig:ex0}.
%
\begin{figure}
   \centering
   \includegraphics[width=6cm]{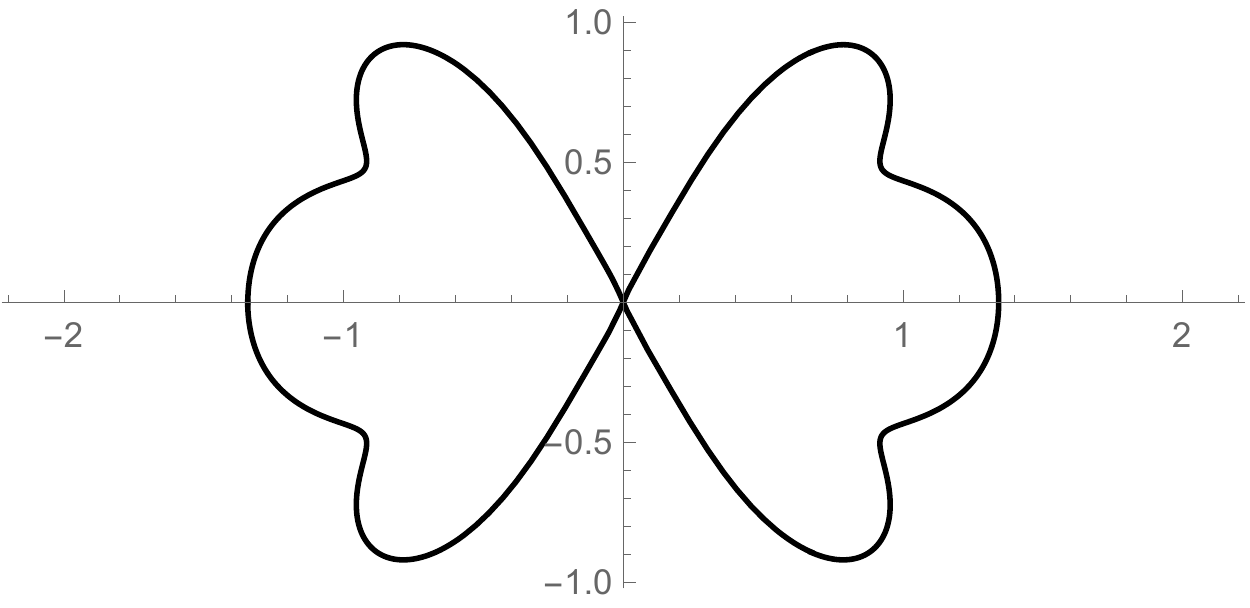}
   \caption{
   Orbit for a solution with $E=0$.
   $(x_0, y_0, v; E)=(0.671188, 0.893818, 0.188131; 6.8 \times 10^{-9})$
    }
   \label{fig:e0}
\end{figure}
We found numerically that this is the only solution with $E=0$. 
The orbit for the $E=0$ solution is shown in figure~\ref{fig:e0} together with the parameters,
which is the member of the series $\gamma$ and resembles the solution shown in figure~\ref{fig:c} (c).
%

Thus around $E=0$, 
there are six solutions for $E>0$ and one solution for $E=0$. 
For $E<0$ there will be many, probably infinitely many, solutions though only five solutions are shown in figure~\ref{fig:ex0}.

The range of the total energy for the figure-eight solutions is
\begin{equation}
  -\frac{5546}{10924}
   \le E \le 0.295542.
\end{equation}
The upper limit is found in three series of solutions $\alpha$, $\beta$ and $\gamma$ in $E \ge 0$,
and is given by a solution
$(x_0, y_0, v; E)$=$(0.686512, 0.639267, 0.646723; 0.295542)$
in the series of solution $\alpha$.
The lower limit is 
the minimum of the potential energy $U$ either for the isosceles triangle configuration 
\begin{equation}
  \min_{x_0, y_0}(u(2y_0)+2u(\sqrt{9x_0^2+y_0^2}))=3u(r_0) = -\frac{3}{4} = -0.75,
\end{equation}
at $(x_0,y_0)=(2^{-5/6}/\sqrt{3},2^{-5/6})$, or for the Euler configuration
\begin{equation}
  \min_{r}(2u(r)+u(2r)) = -\frac{5546}{10924} = -0.507781
\end{equation}
at $r=(2731/1376)^{1/6}=1.121$.

When $x_0 \to \infty$ 
the total energies $E$'s for all series in figure~\ref{fig:ex0}
seem to go toward zero.
For homogeneous system with potential (\ref{eq:homo}), 
by substituting (\ref{eq:alpha6}) into 
(\ref{eq:isos1}) and (\ref{eq:isos2}),
$E$ for figure-eight solution is given by
\begin{equation}
  E = \frac{
  0.0467827
  }{x_0^6}.
\label{eq:alpha6e}
\end{equation}
For the series of solutions $\alpha$,
from the plot of $E x_0^6$ against $x_0$, we can see $E x_0^6$ tends to constant, and  
we obtain asymptotic forms 
for both branches
as 
\begin{equation}
  E \to \frac{0.047}{x_0^6} \; \mbox{or} \;  \frac{0.035}{x_0^6},
\label{eq:LJae}
\end{equation}
by fitting at $x_0=2$, $2.05$ and $2.1$. 
For the other series since we do not have numerical solutions for $x_0 > 1.2$ yet, 
we can not investigate the asymptotic behaviors of $E$ precisely.

\begin{figure}
   \centering
   \includegraphics[width=8cm]{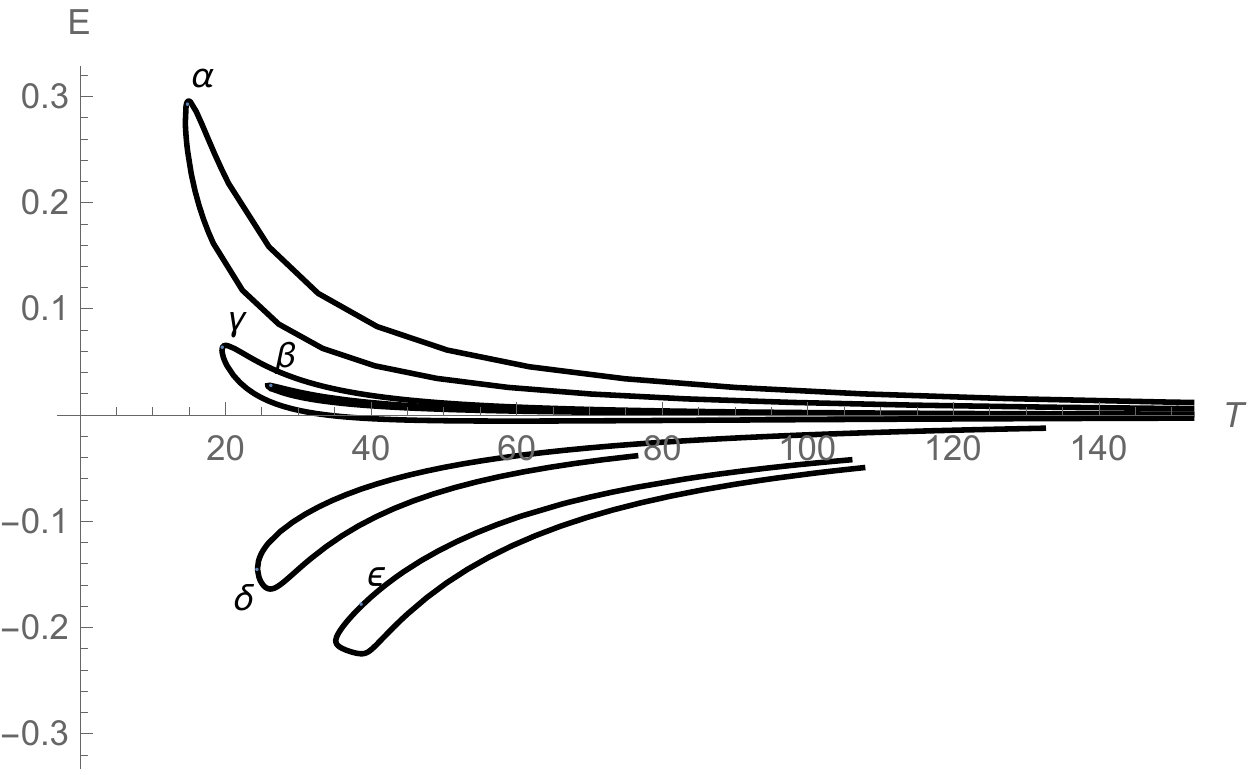} 
   \caption{
  The set $(T,E)$ of the solutions $\alpha$ -- $\epsilon$ in the $(T,E)$-plane.
   }
   \label{fig:TE}
\end{figure}

In figure~\ref{fig:TE}, 
the total energy $E$ for all series of the solutions $\alpha$ -- $\epsilon$ 
are plotted against the period $T$ instead of $x_0$ in $(T,E)$-plane.
We expect a minimum period $T$ within a series of solutions will be smaller for a series with simpler shape of orbit.
Actually in figure~\ref{fig:TE}, we can see a trend that the series with smaller collision number
which may be index of complexity of the shape of the orbit, have smaller minimum period.
Because the series $\alpha$ has the simplest shape of the orbit, 
we conjecture that minimum of $T$ among all solutions 
is the minimum of $T$ in the series $\alpha$, then we find 
\begin{equation}
  14.5 \le T.
\label{eq:Tmin}
\end{equation}

\section{Summary and discussions}
\label{sec:summary}
We have studied on figure-eight choreographic solutions to 
a system of three identical particles interacting through a potential of Lennard-Jones-type (\ref{eq:LJ}). 
The lobe of the figure-eight solutions can be complex shapes 
though it needs to be convex 
in the Newtonian three-body problem \cite{fujiwara}.

By numerical search, we found there are a multitude of such solutions.
A series of solution $\alpha$ tend to a figure-eight solution to a system with homogeneous potential (\ref{eq:homo}) 
when $x_0 \to \infty$. 
The rest are very different from it and have several points 
with large curvatures, collisional segments, in their figure-eight orbits.
These results coincide with theorem by Sbano and Southall \cite{sbano} introduced in the section \ref{sec:intro}.
According to their theorem, there exists a lower bound in the period $T$ of the figure-eight solution, 
which we conjectured in equation (\ref{eq:Tmin}).

In this paper, we have defined the figure-eight choreography as a motion 
starting with the isosceles triangle configurations 
(\ref{eq:isos1}) and (\ref{eq:isos2}),
and going to the Euler configuration 
(\ref{eq:euler1}) and (\ref{eq:euler2}).
%
%
In numerical calculation, however we stop integration in a collinear configuration, 
which means that any collinear configuration is not allowed 
other than the Euler configuration  
(\ref{eq:euler1}) and (\ref{eq:euler2}).
This artificial limitation is introduced by simplicity of our numerical calculation in the first time,
and this should be removed in future. 

In figure~\ref{fig:contour}, we can see very fine structure of the solid curves, $P(x_0,y_0,v)=0$, 
in the $v<0.2$ region. 
\begin{figure}
   \centering
   \includegraphics[width=6cm]{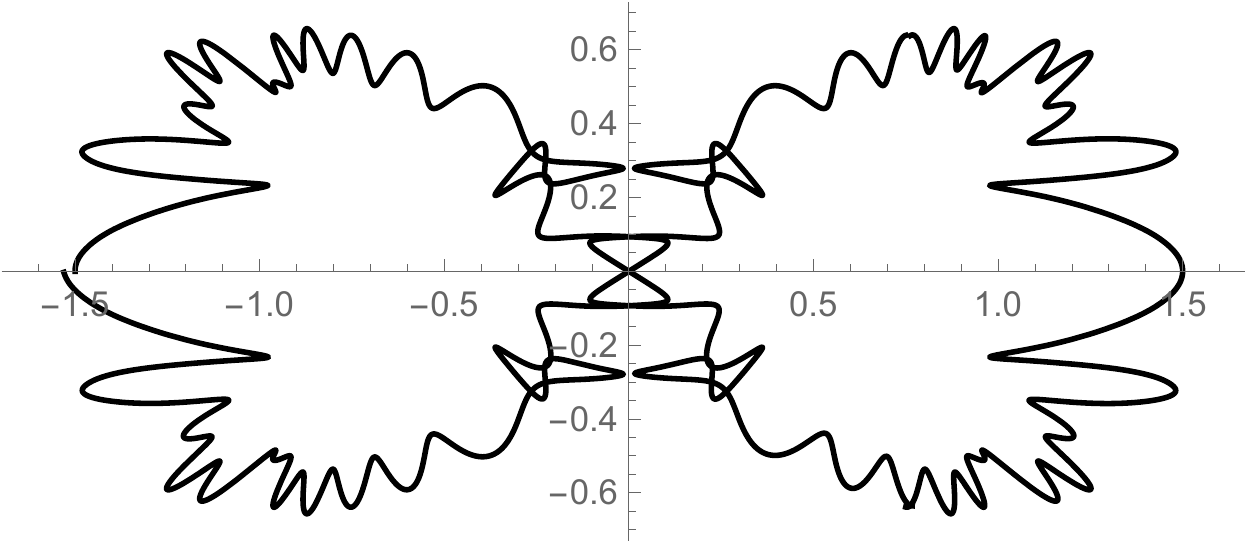}
   \caption{
   An orbit with complex shape for a solution labeled by $\zeta$ in figure~\ref{fig:contour}.
   $(x_0, y_0, v; E)=(0.75, 0.638710, 0.0817816; -0.181634)$
   }
   \label{fig:complex}
\end{figure}
If we magnify this region we see a lot of crossing points of solid and dashed curves,
which mean there are a lot of figure-eight solutions.
An orbit for one of these points labeled by $\zeta$ in figure~\ref{fig:contour} is shown in figure~\ref{fig:complex}.
At $x=-1.5$ curve is not continuous 
because of insufficient accuracy of numerical integration.
In order to explore these complex solutions more elaborate numerical calculations are necessary.
We expect there are infinitely many solutions 
outside the region (\ref{eq:strong}).
Since the solutions originate from the repulsions between bodies 
their orbits must include two points whose distance is less than $r_0$.

Though in the usual shape of figure-eight such as shown in figure~\ref{fig:alpha6} and all orbits shown in this paper
except for figure~\ref{fig:complex}, $\max_{i,t}(x_i(t))=2x_0$ is satisfied, 
this is not true in general. 
We can imagine a counter example from the orbit shown in figure~\ref{fig:complex}.

The importance of the three-body problem under Lennard-Jones-type potential in physical chemistry
and of choreographies eludes us at present.
However, choreographies found in this paper 
may play some role in the molecular dynamics. 
In molecular dynamics\cite{MD},
motions of molecules or atoms are calculated in classical mechanics sometimes 
under Lennard-Jones-type potential, 
that is, $N$-body problem under Lennard-Jones-type potential.
Thus we notice that in such $N$-body problem 
some three bodies could form bound states in positive energies 
such as choreographies with positive energies in $\alpha$, $\beta$ and part of $\gamma$ series we found,
which might have some effect in molecular dynamics. 

Investigation of the figure-eight choreographies under 
various other inhomogeneous potential, such as 
Lennard-Jones-type potential with different powers 
\begin{equation*}
  1/r^b -1/r^a, \; (b,a)\ne(12,6), \; b>a, 
\end{equation*}
Buckingham potential
\begin{equation*}
  e^{-r}-1/r^6,
\end{equation*}
Morse potential
\begin{equation*}
  (1-e^{-a(r-r_0)})^2,
\end{equation*}
or screened Coulomb potential
\begin{equation*}
-e^{-ar}/r,
\end{equation*}
etc., are also interesting.
Does the qualitative nature of our numerical results change dramatically? 
For example, do we still see just one figure-eight with convex lobes?
Do we get an infinite number coalescing with more and more kinks?
These questions should be studied in future, and 
our numerical method exploring the figure-eight solutions 
explained in section \ref{sec:inhomo} 
will be applied for such potentials as an effective numerical method.

Our numerical method has two good features,
compared to the other numerical methods used to find figure-eight orbits, 
such as Moore's relaxation method \cite{moore} or 
truncated Fourier series method by Simo \cite{simo}.
Firstly,
our method can find solution with non minimal of action functional.
This is important since
Sbano and Southall proved that at least one figure-eight solution 
is a mountain-pass critical point of the action functional 
for sufficiently large period $T$ in Lenard-Jones-type potentials \cite{sbano}.
Secondly, in our method, solutions are visualized by contour map of $P(x_0, y_0, v)=0$ and $D(x_0, y_0, v)=0$ 
such as figure~\ref{fig:contour}. 
This is really helpful to explore the solutions
and to understand the relations between solutions. 



\Bibliography{9}


\bibitem{moore}
Moore~C, 1993
{\it Braids in Classical Gravity},
Phys. Rev. Lett. {\bf 70}, 3675--3679

\bibitem{chenAndMont}
Chenciner A and Montgomery R 2000
{\it A remarkable periodic solution of the three-body problem in the case of equal masses},
{\it Annals of Mathematics} {\bf 152}, 881--901

\bibitem{sbano2005}
Sbano L 2005
{\it Symmetric solutions in molecular potentials},
{\it Proceedings of the international conference SPT2004, Symmetry and perturbation theory}, 
(World Scientific Publishing, Singapore) 291--299.

\bibitem{sbano}
Sbano L and Southall J 2010
{\it Periodic solutions of the N-body problem with Lennard-Jones-type potentials},
{\it Dynamical Systems} {\bf 25}, 53--73

\bibitem{lem}
Fujiwara T, Fukuda H and Ozaki H 2003
{\it Choreographic three bodies on the lemniscate},
J. Phys. A: Math. Gen. {\bf 36}, 2791--2800. 

\bibitem{fujiwara}
Fujiwara T and Montgomery R 2005
{\it Convexity in the figure eight solution to the three-body problem},
Journal of Mathematics {\bf 219}, 271--283.

\bibitem{MD} 
Leach A R 
{\it Molecular Modelling},
(Pearson Education; Second Edition 2001)




\bibitem{simo}
Sim\'o C 2001 
{\it Periodic orbits of the planar N-body problem with equal masses and all bodies on the same path},
{\it The Restless Universe: Applications of N-Body Gravitational
Dynamics to Planetary, Stellar and Galactic Systems}, 265-–284, 
ed. Steves B and Maciejewski A, 
NATO Advanced Study Institute, 
(IOP Publishing, Bristol 2001)







\endbib
\end{document}